\documentclass[11pt,reqno,a4paper]{amsart}
\usepackage{amsfonts,amsmath,amssymb}
\usepackage{a4wide}
\usepackage{graphicx}

\allowdisplaybreaks

\newtheorem{theorem}{Theorem}

\newtheorem{coroll}{Corollary}
\theoremstyle{definition}
\newtheorem{remark}{Remark}

\title{Parking function varieties for combinatorial tree models}

\author[A.~Panholzer]{Alois Panholzer}
\address{Alois Panholzer\\
Institut f{\"u}r Diskrete Mathematik und Geometrie\\
Technische Universit\"at Wien\\
Wiedner Hauptstr. 8-10/104\\
1040 Wien, Austria} \email{Alois.Panholzer@tuwien.ac.at}

\date{\today}


\begin{document}

\begin{abstract}
We study the enumeration problem for different kind of tree parking functions introduced recently, called tree parking functions, tree parking distributions, prime tree parking functions, and prime tree parking distributions, for rooted labelled trees of important combinatorial tree families including labelled ordered, unordered and binary trees. Using combinatorial decompositions of the underlying structures yields, after solving the resulting equations, implicit characterizations of suitable generating functions of the total number of such tree parking functions for trees of size $n$ and $n$ successful drivers, from which we obtain exact and asymptotic enumeration results. The approach can be extended to the general situation of tree parking functions for trees of size $n$ and $m<n$ drivers for which we are also able to characterize the generating functions solutions, which allow, by applying analytic combinatorics techniques, a study of the asymptotic behaviour of the total number of tree parking functions and distributions for $n \to \infty$ depending on the load factor $0 < \alpha = \frac{m}{n} < 1$.
\end{abstract}

\maketitle

\section{Introduction}

Parking functions have been introduced by Konheim and Weiss~\cite{KonWei1966} in connection with the analysis of linear probing hashing schemes and since then have become extensively studied combinatorial objects with connections to various other structures as forests, hyperplane arrangements, acyclic functions and non-crossing partitions, see the survey of Yan~\cite{Yan2015} and references therein. Following the vivid description of \cite{KonWei1966} we may think of $n$ free parking spaces $1, 2, \dots, n$ in a row alongside a one-way street, where $n$ drivers arrive sequentially and the $i$-th driver, for $1 \le i \le n$, wishes to park at its preferred parking space $s_{i} \in [n]$ (with $[n]:=\{1, 2, \dots, n\}$). If the parking space $s_{i}$ is free, the $i$-th driver parks there, otherwise he follows the one-way street and parks at the first free parking space, provided that there is such one; otherwise he leaves the street without parking. A parking function of length $n$ is then a sequence $s := (s_{1}, s_{2}, \dots, s_{n}) \in [n]^{n}$ such that all drivers are able to park. When we rank the sequence of preferred parking spaces in non-decreasing order, let us denote this sequence by $s_{(1)} \le s_{(2)} \le \cdots \le s_{(n)}$, then it is easy to see that $s$ is a parking function if and only if $s_{(i)} \le i$, for all $1 \le i \le n$. As an immediate consequence, the order of the preferred parking spaces does not matter, i.e., each permutation of the entries of a parking function $s$ also yields a parking function. A parking function $s$ with $s_{1} \le s_{2} \le \cdots \le s_{n}$ is called an increasing parking function.
It can be shown by numerous ways (again, see \cite{Yan2015} and references) that the number of parking functions of length $n$ is given by $(n+1)^{n-1}$, whereas increasing parking functions of length $n$ are enumerated by the Catalan numbers $B_{n} := \frac{1}{n+1} \binom{2n}{n}$, see~\cite{Sta1999}.

The notion of parking functions has been generalized in various ways; here we focus on generalizations to labelled rooted trees and so-called mappings (i.e., the functional digraphs of functions $[n] \to [n]$) as introduced in \cite{LackPan2016}. We always consider rooted trees with edges oriented towards the root node, where the vertices of any tree $T$ of size $|T|=n$, i.e., with $n$ nodes, are labelled by distinct integers of the set $[n]$; for simplicity we often identify a node with its label. Given a tree $T$ of size $n$ and a sequence $s = (s_{1}, s_{2}, \dots, s_{n}) \in [n]^{n}$ of length $n$, we can adapt the parking procedure described above: the $n$ drivers arrive sequentially and the $i$-th driver wishes to park at its preferred parking space $s_{i}$; starting with node $s_{i}$ he follows the unique path along the directed edges until the first free parking space (i.e., node) is reached, where he parks, or, if there is no such free parking space along the path from $s_{i}$ to the root of $T$, he leaves without parking. 
If each of the $n$ drivers is able to park then we call the pair $(T,s)$ a tree parking function of size $|(T,s)| := |T| = n$. A characterization of tree parking functions has been given in \cite{KinYan2019,LackPan2016}. For a node $v \in T$, let $T_{v}$ be the subtree of $T$ rooted at $v$. Then $(T,s)$ is a tree parking function if and only if $|T_{v}| \le |\{i : s_{i} \in T_{v}\}|$, for all nodes $v \in T$.
Again it follows that any rearrangement of the elements of the sequence $s$ of a tree parking function $(T,s)$ also gives a tree parking function. Tree parking functions $(T,s)$, where the elements of $s$ are forming a non-decreasing sequence, thus generalizations of increasing parking functions, have been introduced in \cite{ButGraYan2017} called tree parking distributions: a tree parking function $(T,s)$ of size $n$, with $s_{1} \le s_{2} \le \cdots \le s_{n}$, is called a tree parking distribution of size $n$. Alternatively, we may consider them as pairs $(T,s)$, with $T$ a size-$n$ tree and $s=\{s_{1}, \dots, s_{n}\}$ a multiset of size $n$ of elements in $[n]$, such that all drivers are able to park.

\begin{figure}
\begin{center}
\text{\raisebox{1.8cm}{\Large{$T$:}}} \quad \includegraphics[height=4cm]{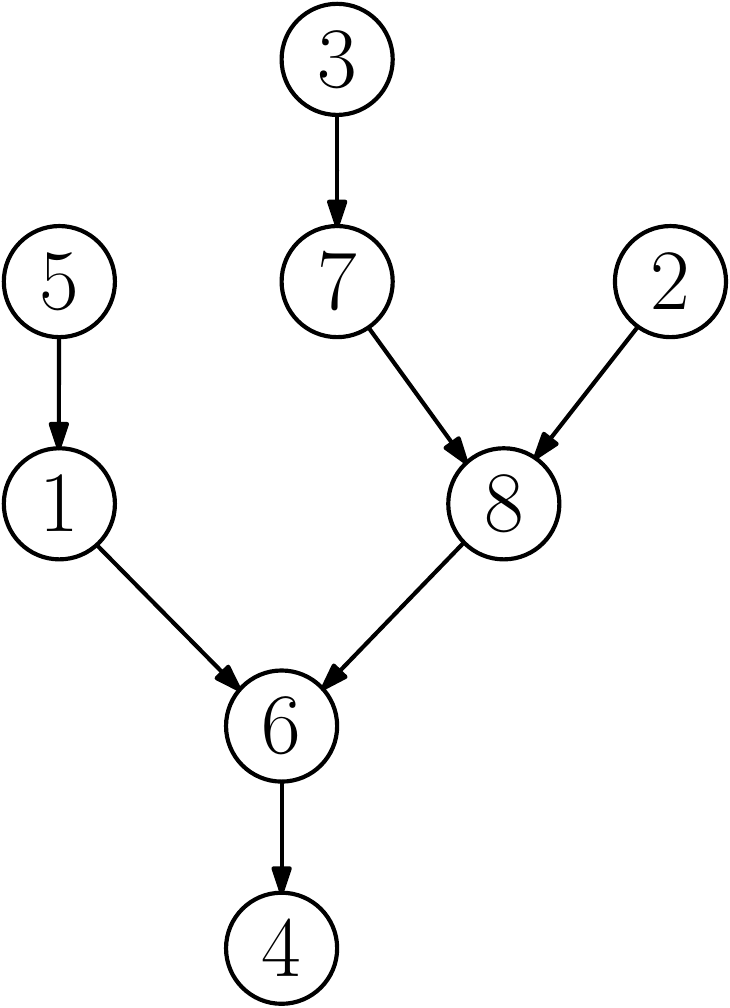}
\caption{A rooted labelled tree $T$ of size $8$. The pair $(T,s)$ with sequence $s=(2,3,5,3,3,1,1,7)$ of preferred parking spaces is a tree parking function, since all drivers are successful. $(T,\tilde{s})$, with $\tilde{s} = (1,1,2,3,3,3,5,7)$, is the corresponding tree parking distribution. During the parking procedure the edges $(5,1)$ and $(2,8)$ are not used, thus this tree parking function and this tree parking distribution are not prime. Let $p=(5,2,3,3,1,2,5,7)$ and $\tilde{p}=(1,2,2,3,3,5,5,7)$, then $(T,p)$ and $(T,\tilde{p})$ are a prime tree parking function and prime tree parking distribution, respectively, since all edges are used during the parking procedure. \label{pic:treeparkingmodels}}
\end{center}
\end{figure}

A further refinement of these notions has been given in \cite{KinYan2019}, where so-called prime tree parking functions and distributions have been introduced and studied. Given a tree parking function $(T,s)$ and a directed edge $e=(v,w)$ in the tree $T$, we say that edge $e$ is used in the parking procedure, if there exists a driver that crosses $e$ during his search for an unoccupied parking space. If all edges in the tree are used in the parking procedure then $(T,s)$ is called a prime tree parking function. Equivalently~\cite{KinYan2019}, a tree parking function $(T,s)$ is prime if it holds $|T_{v}| < |\{i : s_{i} \in T_{v}\}|$, for all non-root nodes $v \in T$. Thus, also each permutation of the elements of the sequence $s$ of a prime tree parking function $(T,s)$ gives a prime tree parking function. Furthermore, prime tree parking functions $(T,s)$, where the elements of the sequence $s$ are forming a non-decreasing sequence, are called prime tree parking distributions.
Figure~\ref{pic:treeparkingmodels} illustrates the different kind of tree parking functions introduced above.

As for classical parking functions, one can extend the notion of tree parking functions and tree parking distributions to the situation that the number $m$ of drivers is less than the number $n$ of parking spaces: a pair $(T,s)$ with $T$ a tree of size $n$ and $s=(s_{1},\dots,s_{m}) \in [n]^{m}$ a sequence of preferred parking spaces is called a $(n,m)$-tree parking function (i.e., tree parking function of size $|T| = n$ and length $|s| = m$) if all $m$ drivers are able to park in the parking procedure. We note that in \cite{LackPan2016} also a characterization of tree parking functions for this extended notion has been given. If the elements of $s$ are forming a non-decreasing sequence (or considering $s$ as a multiset of size $m$ of elements in $[n]$) then we call $(T,s)$ a $(n,m)$-tree parking distribution. For simplicity, we sometimes refer to these objects as general tree parking functions or distributions, respectively.

So far there exist only few enumerative results for tree parking functions and tree parking distributions, and they all concern the family of labelled unordered trees (also called Cayley trees due to the enumeration formula $T_{n} = n^{n-1}$ for the number $T_{n}$ of such trees of size $n$ attributed to A.~Cayley). Namely, in \cite{LackPan2016} the total number $G_{n}$ of tree parking functions of size $n$, i.e., pairs $(T,s)$ with $T$ an unordered labelled tree of size $n$ and $s$ a parking sequence of length $n$, such that all drivers are able to park, has been derived. Moreover, the result has been generalized to the total number $F_{n,m}$ of $(n,m)$-tree parking functions:
\begin{equation*}
\begin{split}
  G_{n} & = \big((n-1)!\big)^{2} \cdot \sum_{k=0}^{n-1} \frac{(n-k)(2n)^{k}}{k!}, \\
  F_{n,m} & = \frac{(n-1)! \, m! \, n^{n-1-m}}{(n-m)!} \cdot \sum_{k=0}^{m} \binom{2m-n-k}{m-k} \frac{(n-k) (2n)^{k}}{k!}.
\end{split}
\end{equation*}
In \cite{KinYan2019}, again for labelled unordered trees, the total number of prime tree parking functions $P_{n}$ and prime tree parking distributions $\tilde{P}_{n}$, respectively, of size $n$ have been computed for which they obtained
\begin{equation*}
  P_{n} = (2n-2)! \qquad \text{and} \qquad  \tilde{P}_{n} = (n-1)! \, S_{n-1},
\end{equation*}
with $S_{n}$ the $n$-th large Schr{\"o}der number.

The aim of this work is to show how the enumeration problem of different notions of tree parking functions and distributions for various combinatorial families $\mathcal{T}$ of rooted labelled trees can be solved in a rather unified way. 
Given a tree family $\mathcal{T}$, let $G_{n}$ and $\tilde{G}_{n}$ be the total number of tree parking functions and tree parking distributions of size $n$, respectively, thus the total number of pairs $(T,s)$, with $T \in \mathcal{T}$ a tree of size $n$ and $s$ a parking sequence or non-decreasing parking sequence, respectively, of length $n$, such that $(T,s)$ is a tree parking function or tree parking distribution, respectively. Furthermore, let $G(z)$ and $\tilde{G}(z)$ be suitably chosen generating functions (g.f.\ for short) of the sequences $(G_{n})_{n}$ and $(\tilde{G}_{n})_{n}$, respectively. Moreover, we denote by $P_{n}$ and $\tilde{P}_{n}$ the total number of prime tree parking functions and prime tree parking distributions, respectively, which are defined in a completely analogous way; suitably chosen generating functions of these sequences are denoted by $P(z)$ and $\tilde{P}(z)$, respectively. 
The starting point of the treatment is a formal recursive description of the family $\mathcal{G}$ of tree parking functions or the family $\mathcal{G}$ of tree parking distributions, respectively, based on a combinatorial decomposition with respect to a certain driver. For tree parking functions we may consider for this decomposition the last driver of the parking sequence; then for several interesting combinatorial tree families including labelled ordered trees, labelled unordered trees and labelled binary trees this yields indeed a suitable recursive description of $\mathcal{G}$. This formal description can be ``translated'' into non-linear differential equations (DEQs) for the generating functions $G(z)$; solving them lead to implicit characterizations of $G(z)$, from which also explicit formul{\ae} for the numbers $G_{n}$ can be obtained. Note that in \cite{LackPan2016} such a decomposition has been applied for labelled unordered trees, but here we avoid the somewhat cumbersome computations based on recurrences for the numbers $G_{n}$, and the use of formal descriptions seem to make the derivations more transparent and easier to extend to further tree families.

For tree parking distributions we extend an idea of \cite{KinYan2019} (applied there for prime tree parking distributions for labelled unordered trees) and decompose the parking distributions w.r.t.\ a driver that arrives at a leaf of the tree; however, unlike in the before-mentioned work, in order to get useful descriptions of the family $\tilde{\mathcal{G}}$ we have to take into account as an auxiliary quantity also the number of leaves in the tree. For various combinatorial tree families, the resulting formal recursive equations for tree parking distributions yield first order non-linear partial differential equations (PDEs) for suitable generating functions $\check{G}(z,v)$, where $v$ encodes the auxiliary quantity. Interestingly, for the tree families considered these PDEs can all be solved by applying the method of characteristics and again yield implicit characterizations of the generating functions as well as explicit results for the numbers $\tilde{G}_{n}$. In particular this allows to enumerate the total number of tree parking distributions for labelled unordered trees.

Furthermore we extend a decomposition of tree parking functions and distributions w.r.t.\ the so-called ``root core'', i.e., the maximal subtree containing the root node, such that all edges contained in the subtree are used during the parking procedure, introduced in \cite{KinYan2019} for labelled unordered trees to more general tree families. This combinatorial decomposition yields relations between the families $\mathcal{G}$ and the corresponding families $\mathcal{P}$ of prime tree parking functions as well as between the families $\tilde{\mathcal{G}}$ and the corresponding families $\tilde{\mathcal{P}}$ of prime tree parking distributions. From these relations we are able to deduce results for the corresponding generating functions $P(z)$ and $\tilde{P}(z)$ and the numbers $P_{n}$ and $\tilde{P}_{n}$, respectively.

Moreover, we show how to obtain results also for the general case, where the number of drivers $m$ is less than the number of parking spaces $n$, thus for the numbers $F_{n,m}$ of $(n,m)$-tree parking functions and the numbers $\tilde{F}_{n,m}$ of $(n,m)$-tree parking distributions. Namely, by decomposing such general tree parking functions and distributions w.r.t.\ either the root node (if it is unoccupied), or the largest subtree containing the root node and only occupied nodes, yields relations between suitable generating functions $F(z,u)$ of $F_{n,m}$ and $\tilde{F}(z,u)$ of $\tilde{F}_{n,m}$, respectively, and the corresponding generating functions $G(z)$ of $G_{n}$ and $\tilde{G}(z)$ of $\tilde{G}_{n}$, respectively. From these relations we are able to deduce implicit characterizations of $F(z,u)$ and $\tilde{F}(z,u)$, from which one could easily gain explicit formul{\ae} for $F_{n,m}$ and $\tilde{F}_{n,m}$. As stated before, in \cite{LackPan2016} the authors obtained for the family of labelled unordered trees a formula for $F_{n,m}$, where they used a different, more involved approach. Furthermore they examined the asymptotic behaviour of $F_{n,m}$, for $n \to \infty$ depending on the ratio $\alpha = \frac{m}{n}$ of the number of drivers and the number of parking spaces (we refer to $\alpha$ as load factor). They observed a phase change behaviour depending on whether $\alpha < \frac{1}{2}$ or $\alpha \ge \frac{1}{2}$, respectively. As a consequence, picking a pair consisting of a tree of size $n$ and a sequence of preferred parking spaces of length $m$, for $\alpha<\frac{1}{2}$, asymptotically, there is a positive probability $p(\alpha) = \frac{\sqrt{1-2\alpha}}{1-\alpha}$ that all drivers are able to park, whereas for $\alpha \ge \frac{1}{2}$, this probability $p(\alpha)$ is $0$. Since then several work \cite{CheGol2019,CurHen2019,GolPrz2019} in the probabilistic literature have observed such phase change phenomena also for generalizations and related problems, and probabilistic explanations for them have been given. Although the focus of the present work is on exact enumeration, we sketch how the generating functions solutions for $F(z,u)$ and $\tilde{F}(z,u)$ can be used to give suitable contour integral representations for $F_{n,m}$ and $\tilde{F}_{n,m}$, respectively, from which one can deduce by an application of the saddle point method also the asymptotic behaviour of the respective number of tree parking functions and tree parking distributions. This yields for each of the tree families considered a phase change behaviour of the probabilities $p(\alpha)$, which is in accordance with very recent general results shown in \cite{CurHen2019}.

\smallskip

The plan of the paper is as follows. In Section~\ref{sec:Results} we introduce the tree families considered in this work and state the main results concerning tree parking functions and distributions, and prime tree parking functions and distributions, which are shown in the subsequent sections. We will always give the proof for the family of ordered labelled trees in more detail, whereas for other tree families we will be more sketchy and omit details.
In Section~\ref{sec:ParkingFunctions} we work out the results for tree parking functions, whereas in Section~\ref{sec:ParkingDistributions} we treat tree parking distributions. Prime tree parking functions and distributions are considered in Section~\ref{sec:Prime}. Generalizations to $(n,m)$-tree parking functions and distributions are given in Section~\ref{sec:General}. Asymptotic results are deduced in Section~\ref{sec:Asymptotic}. We conclude this paper with an outlook and possible extensions in Section~\ref{sec:Outlook}.

\section{Tree models and results\label{sec:Results}}

\subsection{Combinatorial tree families\label{ssec:Tree_families}}

In the following we describe the families $\mathcal{T}$ of labelled trees, for which we enumerate tree parking functions and distributions. For the sake of simplicity we use the term ``combinatorial tree families'', whenever we refer to them throughout this work. As mentioned before, members of these tree families are rooted trees, where each tree of size $n$ is labelled with distinct integers of $[n]$. All of them can be defined in a recursive way.

\begin{description}
\item[Ordered trees] also called planted plane trees. Here the left-to-right ordering of the subtrees of any node in the tree is of relevance. They consist of a root node and a sequence of subtrees, which are, after order-preserving relabellings with labels from $[\text{size of subtree}]$, ordered trees themselves. Thus they are given by the formal equation
\begin{equation*}
  \mathcal{T} = \mathcal{Z} \ast \textsc{Seq}(\mathcal{T}),
\end{equation*}
with $\mathcal{Z}$ an ``atomic element'', i.e., here a vertex, $\ast$ the partition product for labelled combinatorial objects, and $\textsc{Seq}$ the sequence construction (see, e.g., \cite{FlaSed2009} for the combinatorial constructions given here and afterwards).
\item[Unordered trees] also called Cayley trees. Here the left-to-right ordering of the subtrees of any node is not of relevance, thus one could always assume that the subtrees are ordered from left-to-right according to the labels of the root of each subtree. They consist of a root node and a set of subtrees, which are (after suitable order-preserving relabellings) unordered trees themselves. Thus they are given by the formal equation
\begin{equation*}
  \mathcal{T} = \mathcal{Z} \ast \textsc{Set}(\mathcal{T}),
\end{equation*}
with $\textsc{Set}$ the set construction.
\item[$d$-ary trees] Each vertex in the tree has $d$ (distinguishable) positions, where either a node is attached or not. Thus they consist of a root node, and to each of the $d$ positions of the root a $d$-ary tree (after order-preserving relabellings) is attached or not. Formally they can be described via
\begin{equation*}
  \mathcal{T} = \mathcal{Z} \ast \big(\{\epsilon\} + \mathcal{T}\big)^{d},
\end{equation*}
with $\epsilon$ a ``neutral object'' (of size $0$) and $+$ the disjoint union of structures. Sometimes it is convenient to think of $\epsilon$ as an ``external node''. Of particular interest is the case $d=2$, i.e., binary trees, where any node is either a leaf, or there is attached a left child, a right child, or a left and a right child. Furthermore it seems worth to mention that the particular instance $d=1$ corresponds to a ``labelled chain'' of nodes, thus results for parking functions and distributions for $d$-ary trees yield for $d=1$ (apart from a multiplicative factor according to the different labellings of the nodes, i.e., parking spaces) the well-known enumerative results for ordinary parking functions.
\item[$d$-bundled trees] Each vertex in the tree has $d$ (distinguishable) positions, and at each position there can be attached a sequence of nodes. Alternatively, one might think of ordered trees, where for each node the subtrees are separated via $d-1$ walls (or half-edges) into $d$ bundles of subtrees. Thus such trees consist of a root node, and at each of the $d$ positions a sequence of trees is attached, where each of them is (after a suitable relabelling) a $d$-bundled tree itself. We get the formal equation
\begin{equation*}
  \mathcal{T} = \mathcal{Z} \ast \big(\textsc{Seq}(\mathcal{T})\big)^{d}.
\end{equation*}
Of course, since $d=1$ describes ordered trees, $d$-bundled trees are an extension.
\end{description}

\medskip

All of these tree families can be considered as labelled instances of so-called simple families of trees (see \cite{FlaSed2009}), which can be defined via a formal equation of the kind
\begin{equation*}
  \mathcal{T} = \mathcal{Z} \ast \phi\big(\mathcal{T}\big) := \mathcal{Z} \ast \left(\phi_{0} \cdot \{\epsilon\} + \phi_{1} \cdot \mathcal{T} + \phi_{2} \cdot \mathcal{T}^{2} + \cdots\right),
\end{equation*}
with $\phi_{k} \in \mathbb{R}_{0}^{+}$, $\phi_{0} \neq 0$; thus $\phi(\mathcal{T})$ can be considered as a substituted structure reflecting the subtrees of the root node with weights $\phi_{k}$ depending on the degree $k$ of the root.
Alternatively, each tree $T \in \mathcal{T}$ could be regarded as a weighted labelled ordered tree, where the multiplicative weight $w(T)$ of $T$ is given by the product of the weights of its nodes, with the weight of each node $v \in T$ determined by the in-degree\footnote{In combinatorial literature, the edges of rooted trees are considered usually as oriented away from the root node, thus taking the out-degree $d^{+}(v)$ of a node $v$ when introducing simple families of trees, but in the present work the orientation of the edges is reversed.} $d^{-}(v)$, i.e., $w(T) = \prod_{v \in T} \phi_{d^{-}(v)}$. Thus $\phi_{k}$ is the weight of a node with in-degree $k$ and $\phi(t) := \sum_{k \ge 0} \phi_{k} t^{k}$ is called the degree-weight generating function.
The degree-weight generating functions $\phi(t)$ for the tree families considered throughout this work are collected in Table~\ref{tab:combinatorial_tree_families}. Note that for all these tree families the so-called total weight $T_{n} := \sum_{T \in \mathcal{T}, |T|=n} w(T)$ of trees of size $n$ is indeed always an integer, thus $T_{n}$ corresponds to the number of size-$n$ trees of $\mathcal{T}$. These well-known numbers are given in Table~\ref{tab:combinatorial_tree_families}.

In the proposed combinatorial approach enumerating different kinds of tree parking functions for combinatorial tree families, it is important to know what we call here ``the number of attachment points'' $A(T)$ in a tree $T$. Namely, when considering a tree $T \in \mathcal{T}$ of size $|T|=n$, we have to count the number of different trees $T' \in \mathcal{T}$, which can be obtained from $T$ by attaching a new node $v$ (labelled $n+1$) to one of the nodes in $T$. For the tree families introduced before it turns out that the number of attachment points of a tree $T$ depends only on its size $|T|=n$, e.g., it can be seen easily that it is given by $2n-1$ for ordered trees ($v$ can be attached at any node in $T$ left to all other edges or straight to the right of a certain edge), by $n+1$ for binary trees ($v$ can be attached by replacing one of the external nodes in $T$), and by $n$ for unordered trees ($v$ can be attached at one of the nodes in $T$). The number of attachment points of a size-$n$ tree $T$, for the tree-families considered, can be determined easily and are collected in Table~\ref{tab:combinatorial_tree_families}. We remark that for a simple tree family $\mathcal{T}$ with degree-weight generating function $\phi(t)$, the number (or here one might say weight) of attachment points in a tree $T$ is given by the expression $A(T) = \sum_{v \in T} \frac{(d^{-}(v)+1) \, \phi_{d^{-}(v)+1}}{\phi_{d^{-}(v)}}$. In \cite{Pan2006,PanPro2007} (in a different context) a characterization of all simple tree families with the property that $A(T)$ depends only on the size $|T|$ of $T$ has been given. When restricting to tree families, where the total weights $T_{n}$ are always integers, these are exactly the combinatorial tree families stated above.

\begin{table}
\begin{center}
\begin{tabular}{c|c|c|c|c}
& ordered trees & \parbox{1.7cm}{\begin{center}unordered\\[-1mm] trees\end{center}} & $d$-ary trees & $d$-bundled trees\\ 
\hline \hline
degree-weight g.f.\ & & & & \\
$\phi(t)$ & \raisebox{2.5mm}[0mm][0mm]{$\frac{1}{1-t}$} & \raisebox{2.5mm}[0mm][0mm]{$e^{t}$} & \raisebox{2.5mm}[0mm][0mm]{$(1+t)^{d}$} & \raisebox{2.5mm}[0mm][0mm]{$\frac{1}{(1-t)^{d}}$}\\ 
\hline
\# size-$n$ trees & & & & \\
$T_{n}$ & \raisebox{2.5mm}[0mm][0mm]{$(n-1)! \cdot \binom{2n-2}{n-1}$} & \raisebox{2.5mm}[0mm][0mm]{$n^{n-1}$} & \raisebox{2.5mm}[0mm][0mm]{$(n-1)! \cdot \binom{dn}{n-1}$} & \raisebox{2.5mm}[0mm][0mm]{$(n-1)! \cdot \binom{(d+1)n-2}{n-1}$}\\
\hline
\# attachment points & & & & \\
$A(T)$, $|T| = n$ & \raisebox{2.5mm}[0mm][0mm]{$2n-1$} & \raisebox{2.5mm}[0mm][0mm]{$n$} & \raisebox{2.5mm}[0mm][0mm]{$(d-1)n+1$} & \raisebox{2.5mm}[0mm][0mm]{$(d+1)n-1$}\\
\hline
\end{tabular}

\smallskip

\caption{Combinatorial tree families: degree-weight generating functions $\phi(t)$, the number $T_{n}$ of size-$n$ trees, and the number $A(T)$ of attachment points of size-$n$ trees $T$. \label{tab:combinatorial_tree_families}}
\end{center}
\end{table}

Although we mainly focus on stating results for tree parking functions and distributions for these combinatorial tree families, it is possible to extend the approach to further tree families; we comment on this in Section~\ref{sec:ParkingDistributions}.

\subsection{Results\label{ssec:Results}}

The main result of this work concerns implicit characterizations via certain functional equations of suitably defined generating functions of the total number of tree parking functions (of different kinds) for combinatorial tree families.

In the following let us denote by $G(z) := \sum_{n \ge 1} G_{n} \frac{z^{n}}{(n!)^{2}}$ and $P(z) := \sum_{n \ge 1} P_{n} \frac{z^{n}}{(n!)^{2}}$ the double-exponential generating functions of the total number of tree parking functions $G_{n}$ and prime tree parking functions $P_{n}$ of size $n$, respectively, for combinatorial tree families $\mathcal{T}$. Furthermore, let us denote by $\tilde{G}(z) := \sum_{n \ge 1} \tilde{G}_{n} \frac{z^{n}}{n!}$ and $\tilde{P}(z) := \sum_{n \ge 1} \tilde{P}_{n} \frac{z^{n}}{n!}$ the exponential generating functions of the total number of tree parking distributions $\tilde{G}_{n}$ and prime tree parking distributions $\tilde{P}_{n}$ of size $n$, respectively. As mentioned earlier, results for unordered trees are (in parts) already obtained in \cite{KinYan2019,LackPan2016}.

\begin{theorem}\label{thm:gf_treeparking}
Let $X := X(z)$ be one of the generating functions $G := G(z)$, $P := P(z)$, $\tilde{G} := \tilde{G}(z)$, or $\tilde{P} := \tilde{P}(z)$. Then each $X$ can be expressed by means of an auxiliary function $Q := Q(z)$, which is characterized as solution of a certain functional equation of the following form, with some functions $f$ and $\varphi$ depending on the particular tree family\footnote{Note that $\varphi$ characterizing the auxiliary function $Q$ and $\phi$ characterizing the tree family $\mathcal{T}$ do not seem to be related in an obvious way.} and the kind of parking function,
\begin{equation*}
  X = g(Q), \qquad Q = z \varphi(Q),
\end{equation*}
which are specified in Table~\ref{tab:gf_treeparking}.

\begin{table}
\renewcommand{\arraystretch}{1.5}
\begin{equation*}
\begin{array}{c|c|c}
\text{tree family} & \text{generating function} & \text{auxiliary g.f.\ $Q=Q(z)$}\\
\hline \hline
\text{ordered trees} & G = 1-\frac{1}{(1-Q)(1+Q)e^{Q}} & Q = z (1+Q)^{2} e^{Q}\\
\hline
\text{unordered trees \cite{LackPan2016}} & G = 2Q + \ln(1-Q) & Q = z e^{2Q}\\
\hline
\text{$d$-ary trees} & G = \frac{(1-dQ)^{\frac{1}{d}} e^{Q}}{1-Q} - 1 & Q = z \frac{e^{dQ}}{(1-Q)^{d-1}}\\
\hline
\text{$d$-bundled trees} & G = 1-\frac{1}{(1-dQ)^{\frac{1}{d}} (1+Q) e^{Q}} & Q = z (1+Q)^{d+1} e^{dQ}\\
\hline \hline
\text{ordered trees} & P = (1+Q)(1-Q) e^{Q}-1 & Q = \frac{z}{(1-Q)^{2} e^{Q}}\\
\hline
\text{unordered trees \cite{KinYan2019}} & P = 2Q + \ln(1-Q) & Q = \frac{z}{1-Q}\\
\hline
\text{$d$-ary trees} & P = 1-\frac{1-Q}{(1-dQ)^{\frac{1}{d}} \, e^{Q}} & Q = z \frac{e^{Q}}{(1-dQ)^{\frac{d-1}{d}}}\\
\hline
\text{$d$-bundled trees} & P = (1-dQ)^{\frac{1}{d}} (1+Q) e^{Q} -1 & Q = \frac{z}{(1-dQ)^{\frac{d+1}{d}} e^{Q}}\\
\hline
\hline
\text{ordered trees} & \tilde{G} = 1-\frac{1}{(1+Q)^{2} (1-Q)} & Q = z (1+Q)^{4}\\
\hline
\text{unordered trees} & \tilde{G} = Q + \ln(1+Q) + \ln(1-Q) & Q = z (1+Q)^{2} e^{Q}\\
\hline
\text{$d$-ary trees} & \tilde{G} = \frac{(1-Q^{2})^{\frac{1}{d}}}{1-\frac{Q}{d}} -1 & Q = \frac{d z (1+Q)^{2}}{(1-\frac{Q}{d})^{d-1}}\\
\hline
\text{$d$-bundled trees} & \tilde{G} = 1-\frac{1}{(1+\frac{Q}{d}) (1-Q^{2})^{\frac{1}{d}}} & Q = d z (1+Q)^{2} (1+\frac{Q}{d})^{d+1}\\
\hline
\hline
\text{ordered trees} & \tilde{P} = (1+Q)^{2} (1-Q)-1 & Q = \frac{z}{(1-Q)^{2}}\\
\hline
\text{unordered trees \cite{KinYan2019}} & \tilde{P} = Q + \ln(1+Q) + \ln(1-Q) & Q = \frac{z (1+Q)}{1-Q}\\
\hline
\text{$d$-ary trees} & \tilde{P} = 1-\frac{1-\frac{Q}{d}}{(1-Q^{2})^{\frac{1}{d}}} & Q = \frac{dz (1+Q)^{\frac{d+1}{d}}}{(1-Q)^{\frac{d-1}{d}}}\\
\hline
\text{$d$-bundled trees} & \tilde{P} = (1+\frac{Q}{d}) (1-Q^{2})^{\frac{1}{d}} -1 & Q = \frac{dz (1+Q)^{\frac{d-1}{d}}}{(1-Q)^{\frac{d+1}{d}}}\\
\hline
\end{array}
\end{equation*}
\caption{Generating functions solutions $G = G(z)$, $P = P(z)$, $\tilde{G} = \tilde{G}(z)$, and $\tilde{P} = \tilde{P}(z)$ of different kinds of tree parking functions for combinatorial tree families.\label{tab:gf_treeparking}}
\end{table}
\end{theorem}

The asymptotic behaviour of the numbers $G_{n}$, $P_{n}$, $\tilde{G}_{n}$ and $\tilde{P}_{n}$ is characterized in the following theorem.

\begin{theorem}\label{thm:asymptotic_treeparking}
  Let $x_{n}$ be given by one of the normalized enumeration sequences of the different kinds of parking functions, i.e., $\frac{G_{n}}{(n!)^{2}}$, $\frac{P_{n}}{(n!)^{2}}$, $\frac{\tilde{G}_{n}}{n!}$, or $\frac{\tilde{P}_{n}}{n!}$. Then $x_{n}$ is, for $n \to \infty$,  asymptotically given by
\begin{equation*}
	  x_{n} \sim C \, n^{-\frac{5}{2}} \rho^{-n},
\end{equation*}
where $\rho$ is the radius of convergence of the corresponding g.f.\ $X$ determined by $\rho = \frac{\tau}{\varphi(\tau)}$, with $\tau$ the smallest positive real solution of the equation $\varphi(\tau) = \tau \varphi'(\tau)$ (with $X$ and $\varphi$ stated in Theorem~\ref{thm:gf_treeparking}), and $C$ some computable constant, which are summarized in Table~\ref{tab:asymptotic_treeparking}. 

\begin{table}
\renewcommand{\arraystretch}{1.5}
\begin{equation*}
\begin{array}{c|c|c|c}
\hline
& \multicolumn{3}{|c}{ \frac{G_{n}}{(n!)^{2}} \sim C \cdot n^{-\frac{5}{2}} \rho^{-n}}\\
\cline{2-4}
\text{tree family} & \tau & \rho & C\\
\hline
\text{ordered trees} & \sqrt{2}-1 & \frac{\sqrt{2}-1}{2 e^{\sqrt{2}-1}} & \frac{(\sqrt{2}-1)^{\frac{3}{2}} (17+12\sqrt{2})}{8 e^{\sqrt{2}-1} \sqrt{\pi}} \\
\hline
\text{unordered trees} & \frac{1}{2} & \frac{1}{2 e} & \frac{\sqrt{2}}{\sqrt{\pi}}\\
\hline
\text{binary trees} & 1-\frac{1}{\sqrt{2}} & \frac{\sqrt{2}-1}{2 e^{2-\sqrt{2}}} & \frac{e^{1-\frac{1}{\sqrt{2}}} (\sqrt{2}+1)}{2^{\frac{5}{4}} \sqrt{\pi}}\\
\hline
\text{$d$-ary trees} & 1-\sqrt{1-\frac{1}{d}} & \frac{\tau (1-\tau)^{d-1}}{e^{d \tau}} & \frac{\sqrt{\tau} e^{\tau}}{2 \sqrt{d-1} (1-d\tau)^{2-\frac{1}{d}} \sqrt{\pi}}\\
\hline
\text{$d$-bundled trees} & \sqrt{1+\frac{1}{d}}-1 & \frac{\tau}{(1+\tau)^{d+1} e^{d\tau}} & \frac{\sqrt{\tau}}{2 \sqrt{d+1} (1-d\tau)^{2+\frac{1}{d}} e^{\tau} \sqrt{\pi}}\\
\hline
\hline
& \multicolumn{3}{|c}{ \frac{P_{n}}{(n!)^{2}} \sim C \cdot n^{-\frac{5}{2}} \rho^{-n}}\\
\cline{2-4}
\text{tree family} & \tau & \rho & C\\
\hline
\text{ordered trees} & \sqrt{2}-1 & 2 (3-2\sqrt{2}) e^{\sqrt{2}-1} & \frac{(3-2\sqrt{2}) e^{\sqrt{2}-1}}{2 \sqrt{\pi}} \\
\hline
\text{unordered trees} & \frac{1}{2} & \frac{1}{4} & \frac{1}{4 \sqrt{\pi}}\\
\hline
\text{binary trees} & 1-\frac{1}{\sqrt{2}} & \frac{(\sqrt{2}-1)^{\frac{3}{2}}}{\sqrt{2} \, e^{1-\frac{1}{\sqrt{2}}}} & \frac{\sqrt{\sqrt{2}-1}}{2 e^{1-\frac{1}{\sqrt{2}}} \sqrt{\pi}}\\
\hline
\text{$d$-ary trees} & 1-\sqrt{1-\frac{1}{d}} & \frac{\tau (1-d\tau)^{1-\frac{1}{d}}}{e^{\tau}} & \frac{\sqrt{d} \sqrt{2d\tau-1} \, \rho}{2 (1-d\tau)^{2} \sqrt{\pi}}\\
\hline
\text{$d$-bundled trees} & \sqrt{1+\frac{1}{d}}-1 & \tau (1-d\tau)^{1+\frac{1}{d}} e^{\tau} & \frac{\sqrt{d} \sqrt{1-2d\tau} \, \rho}{2 (1-d\tau)^{2} \sqrt{\pi}}\\
\hline
\hline
& \multicolumn{3}{|c}{ \frac{\tilde{G}_{n}}{n!} \sim C \cdot n^{-\frac{5}{2}} \rho^{-n}}\\
\cline{2-4}
\text{tree family} & \tau & \rho & C\\
\hline
\text{ordered trees} & \frac{1}{3} & \frac{27}{256} & \frac{27 \sqrt{6}}{128 \sqrt{\pi}} \\
\hline
\text{unordered trees} & \sqrt{2}-1 & \frac{\sqrt{2}-1}{2 e^{\sqrt{2}-1}} & \frac{\sqrt{\sqrt{2}-1} \, (3+2\sqrt{2})}{4\sqrt{\pi}}\\
\hline
\text{binary trees} & \frac{1}{2} & \frac{1}{12} & \frac{\sqrt{3}}{2 \sqrt{\pi}} \\
\hline
\text{$d$-ary trees} & \frac{\sqrt{2d(d-1)}-2}{d-2}-1 & \frac{\tau(1-\frac{\tau}{d})^{d-1}}{d(1+\tau)^{2}} & \text{\parbox{6cm}{$\frac{\big((2d^{2}-5d+2)\tau+5d^{2}-4d\big) \sqrt{(d^{2}+d+2)\tau-2d}}{4d^{2}(d-1)(d-2)(1-\frac{\tau}{d})^{2} \sqrt{\pi}}$\\ \hspace*{5mm} $\times (1+\tau)^{\frac{1}{d}}(1-\tau)^{\frac{1}{d}}$}} \\
\hline
\text{$d$-bundled trees} & \frac{2+\sqrt{2d(d+1)}}{d+2}-1 & \frac{\tau}{d(1+\tau)^{2} (1+\frac{\tau}{d})^{d+1}} & \frac{\big((2d^{2}+5d+2)\tau+5d^{2}+4d\big)\sqrt{(d^{2}-d+2)\tau+2d}}{4d^{2}(d+1)(d+2)(1+\frac{\tau}{d})^{2} (1+\tau)^{\frac{1}{d}} (1-\tau)^{\frac{1}{d}} \sqrt{\pi}}\\
\hline
\hline
& \multicolumn{3}{|c}{ \frac{\tilde{P}_{n}}{n!} \sim C \cdot n^{-\frac{5}{2}} \rho^{-n}}\\
\cline{2-4}
\text{tree family} & \tau & \rho & C\\
\hline
\text{ordered trees} & \frac{1}{3} & \frac{4}{27} & \frac{2 \sqrt{3}}{27 \sqrt{\pi}} \\
\hline
\text{unordered trees} & \sqrt{2}-1 & 3-2\sqrt{2} & \frac{\sqrt{3\sqrt{2}-4}}{2 \sqrt{\pi}}\\
\hline
\text{binary trees} & \frac{1}{2} & \frac{\sqrt{3}}{18} & \frac{\sqrt{2}}{6 \sqrt{\pi}}\\
\hline
\text{$d$-ary trees} & \frac{\sqrt{2d(d-1)}-2}{d-2}-1 & \frac{\tau (1-\tau)^{1-\frac{1}{d}}}{d (1+\tau)^{1+\frac{1}{d}}} & \frac{\sqrt{d((3d+2)\tau-d-2)} \, \rho}{(d-2)(1-\tau)^{2} \sqrt{\pi}}\\
\hline
\text{$d$-bundled trees} & \frac{2+\sqrt{2d(d+1)}}{d+2}-1 & \frac{\tau (1-\tau)^{1+\frac{1}{d}}}{d(1+\tau)^{1-\frac{1}{d}}} & \frac{\sqrt{d((3d-2)\tau-d+2)} \, \rho}{(d+2)(1-\tau)^{2} \sqrt{\pi}}\\
\hline
\end{array}
\end{equation*}
\caption{Asymptotic behaviour of the counting sequences $G_{n}$, $P_{n}$, $\tilde{G}_{n}$, and $\tilde{P}_{n}$ of different kinds of tree parking functions for combinatorial tree families.\label{tab:asymptotic_treeparking}}
\end{table}
\end{theorem}

\begin{remark}
It is somewhat surprising that the sublinear term in the asymptotic behaviour of the number of tree parking functions (of different kinds) is different from $n^{-\frac{3}{2}}$, which is ``typical'' for the enumeration of tree families, but it is $n^{-\frac{5}{2}}$ and thus the one occurring often, e.g., in the enumeration of so-called maps (see, e.g., \cite{FlaSed2009}).
\end{remark}

From the generating functions solutions given in Theorem~\ref{thm:gf_treeparking} one can easily obtain explicit formul{\ae} for the underlying enumeration sequences. For a few instances one even gets well-known sequences occurring in the online encyclopedia of integer sequences (OEIS)~\cite{OEIS2020}. In the following we give such results, where we restrict ourselves to the most prominent combinatorial tree families, i.e., ordered trees, unordered trees, and binary trees.

\begin{coroll}\label{cor:explicitformulae_treeparking}
The numbers $G_{n}$, $P_{n}$, $\tilde{G}_{n}$ and $\tilde{P}_{n}$ enumerating different kinds of tree parking functions are for  ordered trees, unordered trees, and binary trees given by the explicit enumeration formul{\ae} stated in Table~\ref{tab:explicitformulae_treeparking}.

\begin{table}
\renewcommand{\arraystretch}{1.5}
\begin{equation*}
\begin{array}{c|c|c}
\text{tree family} & \text{enumeration formula} & \text{\parbox{2.5cm}{\begin{center}relation to\\ OEIS sequence\end{center}}}\\
\hline \hline
\text{ordered trees} & \text{\parbox{7cm}{$G_{n} = n! (n-2)! \sum_{k=0}^{n-2} \frac{(n-1)^{k}}{k!}$\\ $\times \sum_{\ell=0}^{n-2-k} (\ell+1)(2\ell+3) \binom{2n-3}{n-2-k-\ell}$}} & \\
\hline
\text{unordered trees \cite{LackPan2016}} & G_{n} = ((n-1)!)^{2} \sum_{k=0}^{n-1} \frac{(n-k)(2n)^{k}}{k!} & \\
\hline
\text{binary trees} & \text{\parbox{9cm}{$G_{n} = n! (n-1)! \sum_{k=0}^{n-1} \frac{(2n+1)^{k}}{k!}$\\ $\times \sum_{\ell=0}^{n-1-k} \frac{1}{2^{\ell}} \binom{2\ell}{\ell} \binom{2n-k-\ell -2}{n-1} \left[2-\frac{(2n-k-\ell)(2n-k-\ell-1)}{n(n+1)}\right]$}} & \\
\hline \hline
\text{ordered trees} & P_{n} = n! (n-2)! \sum_{k=0}^{n-2} \binom{2n-1+k}{k} \frac{(1-n)^{n-2-k}}{(n-2-k)!} & \\
\hline
\text{unordered trees \cite{KinYan2019}} & P_{n} = (2n-2)! & A010050 \\
\hline
\text{binary trees} & P_{n} = n! (n-2)! \sum_{k=0}^{n-2} \binom{\frac{n+1}{2} +k}{k} \frac{2^{k+1} (n-1)^{n-2-k}}{(n-2-k)!} & \\
\hline
\hline
\text{ordered trees} & \tilde{G}_{n} = (n-1)! \sum_{k=0}^{n-1} \binom{4n-3}{k} (3-2(n-k)) & A294084\\
\hline
\text{unordered trees} & \text{\parbox{6cm}{$\tilde{G}_{n} = (n-1)! \sum_{k=0}^{n-1} \frac{n^{k}}{k!}$\\ $\times \left[\frac{3n+k+1}{n+k+1} \binom{2n-1}{n+k} - \sum_{\ell=0}^{n-1-k} \binom{2n}{\ell}\right]$}} & \\
\hline
\text{binary trees} & \text{\parbox{7cm}{$\tilde{G}_{n} = (n-1)! \, 2^{n-1} \sum_{k=0}^{n-1} \binom{2n-\frac{1}{2}}{k}$\\ $\times \sum_{\ell=0}^{n-k-1} \frac{1}{2^{\ell}} \binom{n+\ell+1}{\ell} \binom{n-k-\ell-\frac{3}{2}}{n-k-\ell-1} \frac{1-2(n-k-\ell)}{2(n-k-\ell)-3}$}} & \\
\hline
\hline
\text{ordered trees} & \tilde{P}_{n} = \frac{2(3n-3)!}{(2n-1)!} & A000139\\
\hline
\text{unordered trees \cite{KinYan2019}} & \text{\parbox{8cm}{$\tilde{P}_{n}=2(n-2)! \sum_{k=0}^{n-2} \binom{n+k}{k} \binom{n-1}{k+1} = (n-1)! S_{n-1}$,\\ with $S_{n}$ the large Schr{\"o}der numbers}} & A006318\\
\hline
\text{binary trees} & \tilde{P}_{n} = (n-1)! \frac{2^{2n-1}}{n+1} \binom{\frac{3(n-1)}{2}}{n-1} & A214377\\
\hline
\end{array}
\end{equation*}
\caption{Explicit enumeration formul{\ae} for the numbers $G_{n}$, $P_{n}$, $\tilde{G}_{n}$, and $\tilde{P}_{n}$ of different kinds of tree parking functions for ordered trees, unordered trees, and binary trees.\label{tab:explicitformulae_treeparking}}
\end{table}
\end{coroll}

\begin{remark}
A few enumeration sequences of tree parking functions occurring above we find particularly interesting. First, for the number of prime tree parking distributions on ordered trees it holds $\tilde{P}_{n} = n! A_{n-1}$, with $A_{n}$ given by OEIS $A000139$ enumerating, amongst others, also rooted non-separable planar maps. Furthermore, also for ordered trees it holds that the number of tree parking distributions satisfies $\tilde{G}_{n} = n! A_{n}$, with $A_{n}$ given by OEIS $A294084$ enumerating indecomposable intervals in the Tamari lattices. For prime tree parking distributions on binary trees we get $\tilde{P}_{n} = (n-1)! A_{n-1}$, with numbers $A_{n}$ given by OEIS $A214377$, which occur as coefficients of the g.f.\ $A(z)$ satisfying the equation $A = 1 + 4z A^{\frac{3}{2}}$. It seems that so far there has not been given a concrete combinatorial meaning for the latter numbers and here we could add such one.
\end{remark}

Now we turn to general tree parking functions and distributions and denote by $F(z,u) := \sum_{n \ge 1} \sum_{0 \le m \le n} F_{n,m} \frac{z^{n} u^{n-m}}{n! \, m!}$ and $\tilde{F}(z,u) := \sum_{n \ge 1} \sum_{0 \le m \le n} \tilde{F}_{n,m} \frac{z^{n} u^{n-m}}{n!}$ suitable generating functions of the number of $(n,m)$-tree parking functions $F_{n,m}$ and $(n,m)$-tree parking distributions $\tilde{F}_{n,m}$, respectively, for combinatorial tree families $\mathcal{T}$. As a main result, also for the general case we get implicit characterizations of these generating functions via certain functional equations (extending the corresponding ones given in Theorem~\ref{thm:gf_treeparking}, which can be obtained by setting $u=0$). As mentioned earlier, the result for tree parking functions for the family of unordered trees has been obtained in \cite{LackPan2016}.

\begin{theorem}\label{thm:gf_general_treeparking}
Let $X := X(z,u)$ be one of the generating functions $F := F(z,u)$ or $\tilde{F} := \tilde{F}(z,u)$. Then each $X$ can be expressed by means of an auxiliary function $\hat{Q} := \hat{Q}(z,u)$, which is characterized as solution of a certain functional equation of the following form, with some functions $f$ and $\varphi$ depending on the particular tree family and the kind of parking function,
\begin{equation*}
  X = f(Q,u), \qquad Q = z \varphi(Q,u),
\end{equation*}
which are specified in Table~\ref{tab:gf_general_treeparking}.

\begin{table}
\renewcommand{\arraystretch}{1.5}
\begin{equation*}
\begin{array}{c|c|c}
\text{tree family} & \text{generating function} & \text{auxiliary g.f.\ $\hat{Q}=\hat{Q}(z,u)$}\\
\hline \hline
\text{ordered trees} & F = 1-\frac{1}{(1-\hat{Q})(1+\hat{Q}+u\hat{Q}(1-\hat{Q}))e^{\hat{Q}}} & \hat{Q} = z \big(1+\hat{Q}+u\hat{Q}(1-\hat{Q})\big)^{2} e^{\hat{Q}}\\
\hline
\text{unordered trees} & F = \ln(1-\hat{Q}) + \hat{Q} (2+u(1-\hat{Q})) & \hat{Q} = z e^{\hat{Q} (2+u(1-\hat{Q}))}\\
\hline
\text{$d$-ary trees} & F = \frac{(1-d\hat{Q})^{\frac{1}{d}} e^{\hat{Q}}}{1-\hat{Q}-u\hat{Q}(1-d\hat{Q})}-1 & \hat{Q} = \frac{z \, e^{d \hat{Q}}}{(1-\hat{Q}-u\hat{Q}(1-d\hat{Q}))^{d-1}}\\
\hline
\text{$d$-bundled trees} & F = 1-\frac{1}{(1-d\hat{Q})^{\frac{1}{d}} (1+\hat{Q}+u\hat{Q}(1-d\hat{Q})) e^{\hat{Q}}} & \hat{Q} = z \big(1+\hat{Q}+u\hat{Q}(1-d\hat{Q})\big)^{d+1} e^{d \hat{Q}}\\
\hline
\hline
\text{ordered trees} & \tilde{F} = 1 - \frac{1}{(1-\hat{Q}) ((1+\hat{Q})^{2} + u\hat{Q}(1-\hat{Q}))} & \hat{Q} = z \big((1+\hat{Q})^{2} +u\hat{Q}(1-\hat{Q})\big)^{2}\\
\hline
\text{unordered trees} & \text{\parbox{5.5cm}{$\tilde{F} = \hat{Q} + \ln(1+\hat{Q}) + \ln(1-\hat{Q})$ \\ \hspace*{0.6cm} $\mbox{} + \frac{u \hat{Q} (1-\hat{Q})}{1+\hat{Q}}$}} & \hat{Q} = z (1+\hat{Q})^{2} e^{\hat{Q} \big(1+\frac{u(1-\hat{Q})}{1+\hat{Q}}\big)}\\
\hline
\text{$d$-ary trees} & \tilde{F} = \frac{(1+\hat{Q})^{\frac{1}{d}} (1-\hat{Q})^{\frac{1}{d}}}{1-\frac{\hat{Q}}{d} - \frac{u\hat{Q}(1-\hat{Q})}{d(1+\hat{Q})}} - 1 & \hat{Q} = \frac{d \, z (1+\hat{Q})^{2}}{\big(1-\frac{\hat{Q}}{d} - \frac{u\hat{Q}(1-\hat{Q})}{d(1+\hat{Q})}\big)^{d-1}}\\
\hline
\text{$d$-bundled trees} & \tilde{F} = 1-\frac{1}{(1+\hat{Q})^{\frac{1}{d}} (1-\hat{Q})^{\frac{1}{d}} \big(1+\frac{\hat{Q}}{d} + \frac{u\hat{Q}(1-\hat{Q})}{d(1+\hat{Q})}\big)} & \text{\parbox{5cm}{$\hat{Q} = d \, z (1+\hat{Q})^{2}$ \\ \hspace*{0.6cm} $\times \big(1+\frac{\hat{Q}}{d} + \frac{u\hat{Q}(1-\hat{Q})}{d(1+\hat{Q})}\big)^{d+1}$}}\\
\hline
\end{array}
\end{equation*}
\caption{Generating functions solutions of general tree parking functions $F=F(z,u)$ and distributions $\tilde{F}=\tilde{F}(z,u)$ for combinatorial tree families.\label{tab:gf_general_treeparking}}
\end{table}
\end{theorem}

As mentioned earlier, the numbers $F_{n,m}$ and $\tilde{F}_{n,m}$ of $(n,m)$-tree parking functions and distributions, respectively, change its asymptotic behaviour, for $n, m  \to \infty$, depending on the ratio $\alpha = \frac{m}{n}$. We state these results, where we restrict ourselves to the most important members of combinatorial tree families, ordered, unordered and binary trees.

\begin{theorem}\label{thm:asymptotic_general_treeparking}
  Let $p_{n,m} := \frac{F_{n,m}}{T_{n} \, n^{m}}$ or $p_{n,m} := \frac{\tilde{F}_{n,m}}{T_{n} \, \binom{n+m-1}{m}}$ be the probability that a randomly chosen sequence of length $m$ respectively multiset of size $m$ on the set $[n]$ is, for a randomly chosen size-$n$ tree of one of the families of ordered, unordered, or binary trees, a parking function or parking distribution, respectively. Furthermore, let $p(\alpha) := \lim_{n \to \infty} p_{n, \alpha n}$, for a load factor $0 \le \alpha \le 1$. Then the probabilities $p(\alpha)$ undergo a phase change behaviour of the form
\begin{equation*}
  p(\alpha) = \begin{cases}
	  C(\alpha) \cdot \sqrt{1-\frac{\alpha}{\alpha_{0}}}, & \quad \text{for \: $0 \le \alpha < \alpha_{0}$},\\
		0, & \quad \text{for \: $\alpha_{0} \le \alpha \le 1$},
	\end{cases}
\end{equation*}
with critical load factor $\alpha_{0}$ and function $C(\alpha)$ bounded on $[0,\alpha_{0}]$ that depend on the particular tree family. These results are summarized in Table~\ref{tab:asymptotic_general_treeparking}.

\begin{table}
\begin{equation*}
\begin{array}{c|c|c|c}
\hline
& \multicolumn{3}{|c}{\text{tree parking functions}}\\
\cline{2-4}
\raisebox{2mm}[-2mm]{\text{tree family}} & p(\alpha) & C(\alpha) & \alpha_{0}\\
\hline
\text{ordered trees} & \begin{cases} \frac{\sqrt{1-2\alpha-\alpha^{2}}}{(1-\alpha)^{2} e^{\alpha}}, & 0 \le \alpha < \alpha_{0},\\
0, & \alpha_{0} \le \alpha \le 1 \end{cases} & \frac{\sqrt{(\sqrt{2}-1)\alpha+1}}{(1-\alpha)^{2} e^{\alpha}} & \sqrt{2} -1\\
\hline
\text{unordered trees} & \begin{cases} \frac{\sqrt{1-2\alpha}}{1-\alpha}, & 0 \le \alpha < \alpha_{0},\\
0, & \alpha_{0} \le \alpha \le 1 \end{cases} & \frac{1}{1-\alpha} & \frac{1}{2}\\
\hline
\text{binary trees} & \begin{cases} \frac{\sqrt{2-4\alpha+\alpha^{2}} \, e^{\frac{\alpha}{2}}}{\sqrt{2-2\alpha}}, & 0 \le \alpha < \alpha_{0},\\
0, & \alpha_{0} \le \alpha \le 1 \end{cases} & \frac{\sqrt{2-(2-\sqrt{2})\alpha} \, e^{\frac{\alpha}{2}}}{\sqrt{2-2\alpha}} & 2-\sqrt{2}\\
\hline
\hline
& \multicolumn{3}{|c}{\text{tree parking distributions}}\\
\cline{2-4}
\raisebox{2mm}[-2mm]{\text{tree family}} & p(\alpha) & C(\alpha) & \alpha_{0}\\
\hline
\text{ordered trees} & \begin{cases} \frac{\sqrt{1-3\alpha}}{(1-\alpha)^{2} \sqrt{1+\alpha}}, & 0 \le \alpha < \alpha_{0},\\
0, & \alpha_{0} \le \alpha \le 1 \end{cases} & \frac{1}{(1-\alpha)^{2} \sqrt{1+\alpha}} & \frac{1}{3}\\
\hline
\text{unordered trees} & \begin{cases} \frac{\sqrt{1-2\alpha-\alpha^{2}}}{1-\alpha}, & 0 \le \alpha < \alpha_{0},\\
0, & \alpha_{0} \le \alpha \le 1 \end{cases} & \frac{\sqrt{(\sqrt{2}-1)\alpha+1}}{1-\alpha} & \sqrt{2}-1\\
\hline
\text{binary trees} & \begin{cases} \frac{\sqrt{1-\alpha-2\alpha^{2}}}{\sqrt{1-\alpha}}, & 0 \le \alpha < \alpha_{0},\\
0, & \alpha_{0} \le \alpha \le 1 \end{cases} & \frac{\sqrt{1+\alpha}}{\sqrt{1-\alpha}} & \frac{1}{2}\\
\hline
\end{array}
\end{equation*}
\caption{Phase change behaviour of the limiting probabilities $p(\alpha) = \lim_{n \to \infty} p_{n, \alpha n}$ depending on whether the load factor $\alpha$ is below or above the critical load factor $\alpha_{0}$.\label{tab:asymptotic_general_treeparking}}
\end{table}
\end{theorem}

\section{Parking functions\label{sec:ParkingFunctions}}

\subsection{Combinatorial decomposition\label{ssec:ParkingFunction_Decomposition}}

\begin{figure}
\begin{center}
\raisebox{0.45cm}{\includegraphics[width=3.2cm]{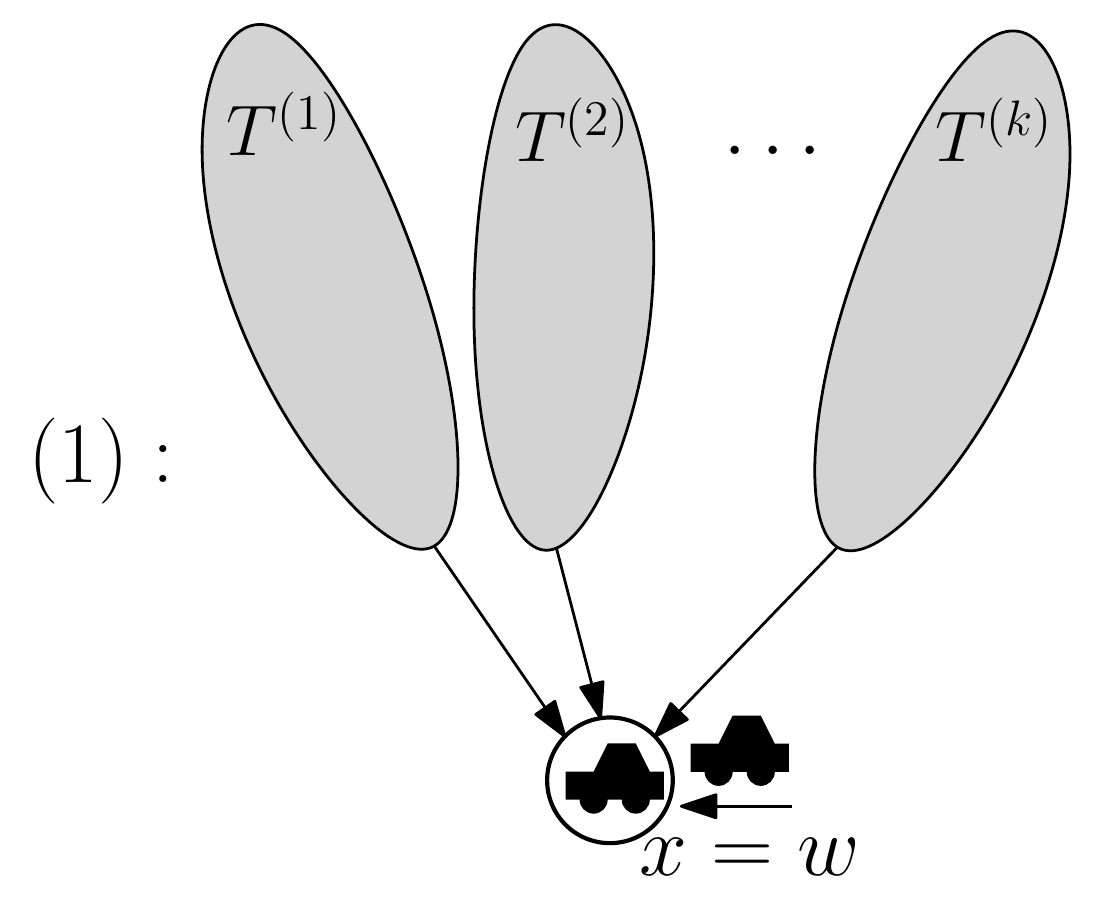}} \;
\raisebox{0.45cm}{\includegraphics[width=3.2cm]{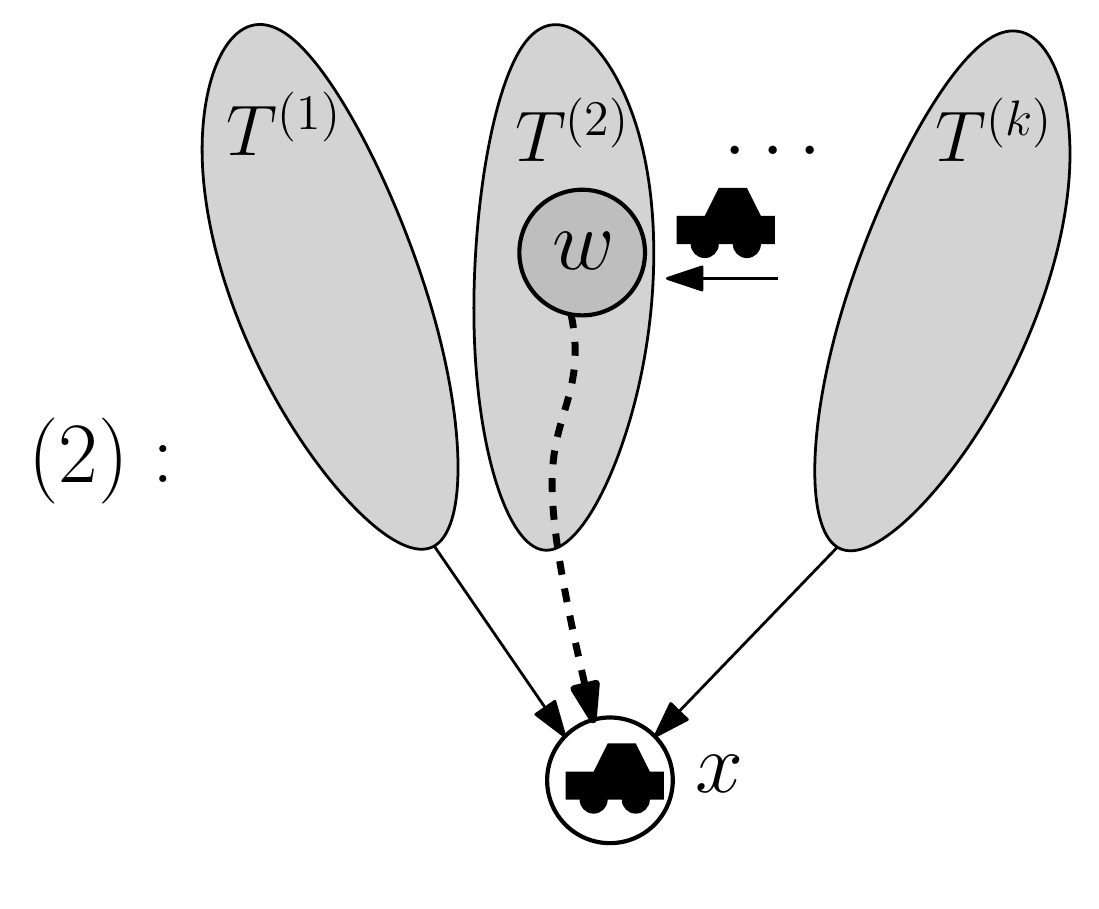}} \;
\includegraphics[width=4cm]{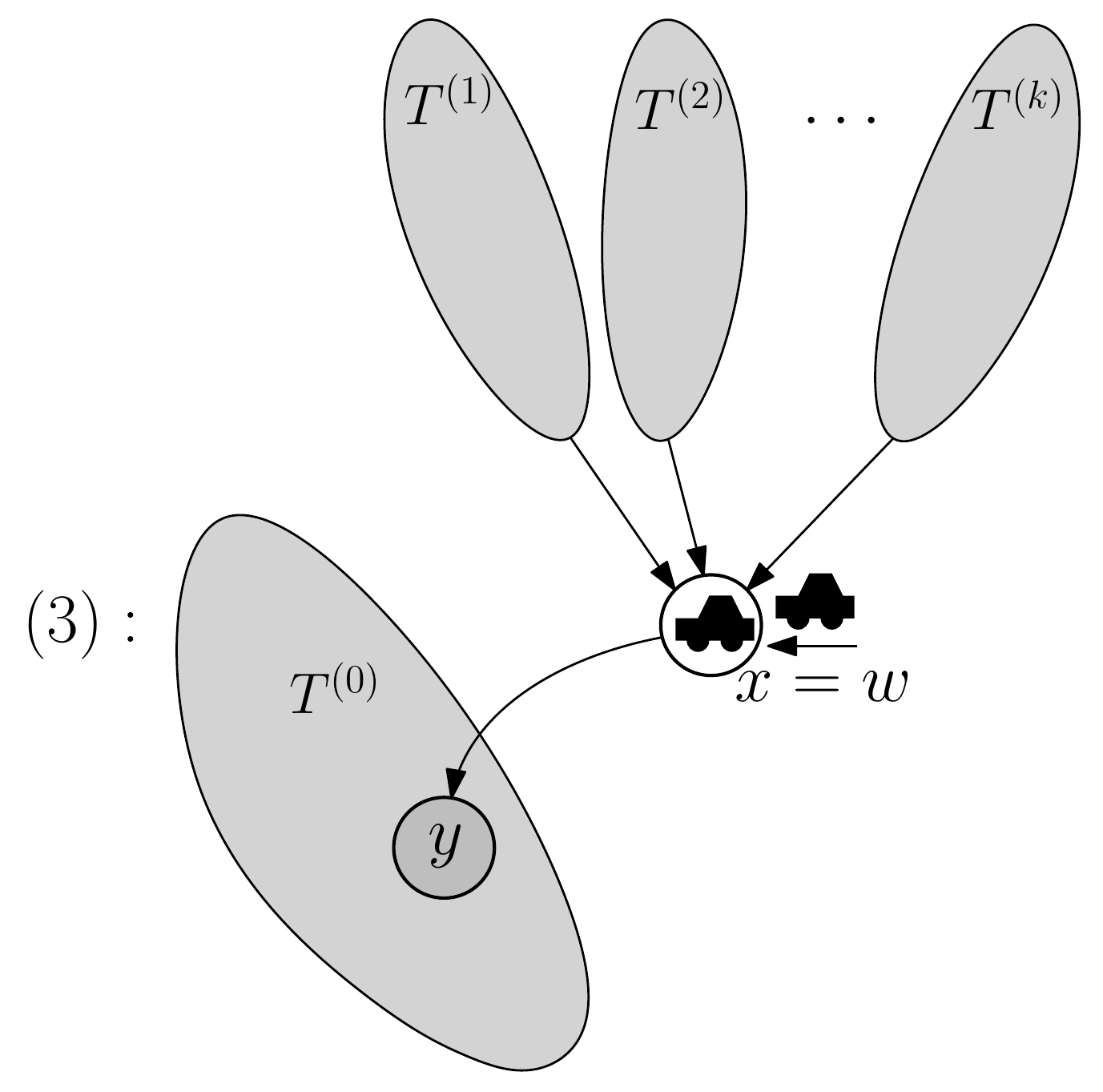} \;
\includegraphics[width=4cm]{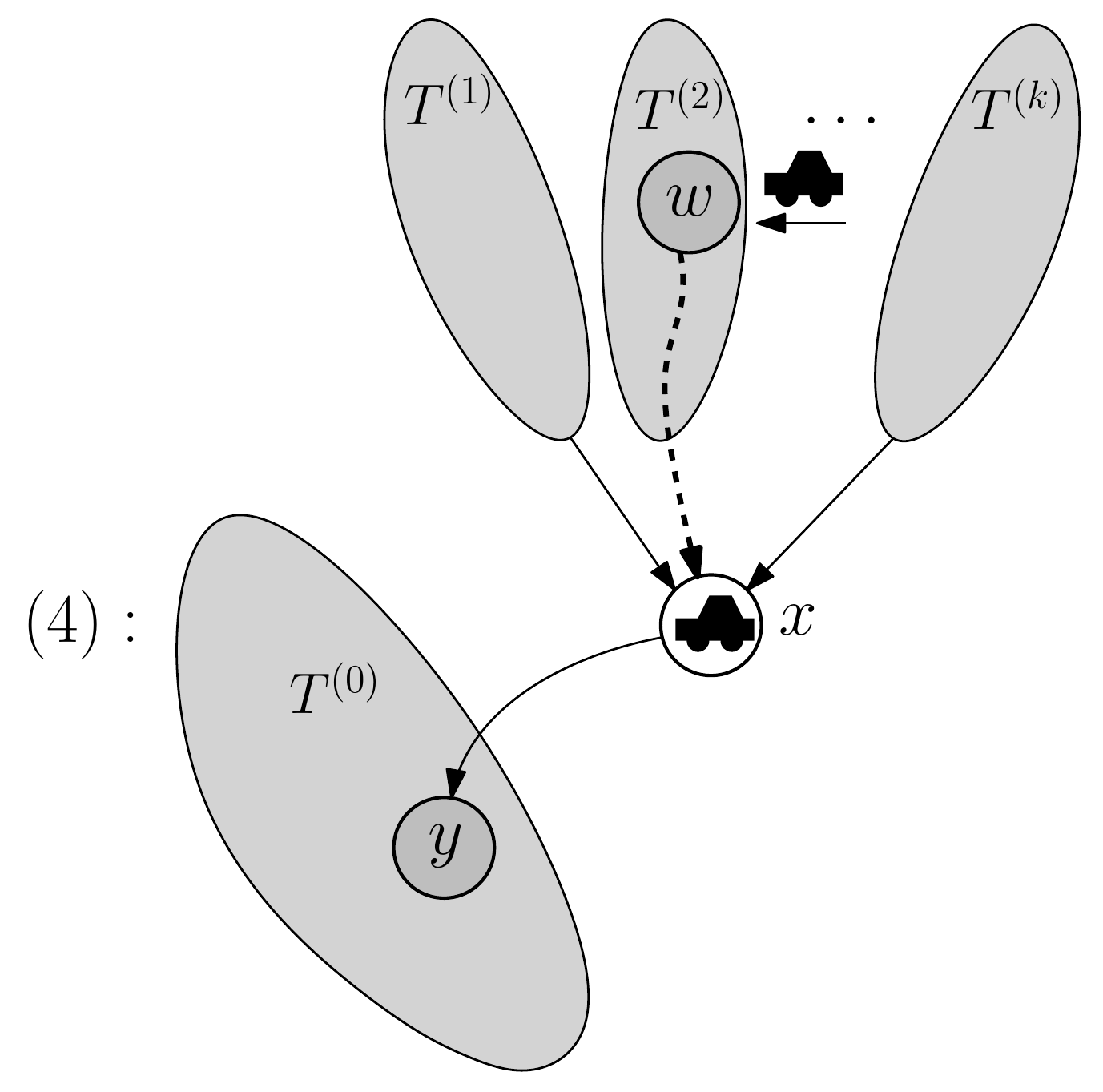}
\caption{Decomposition of a tree parking function $(T,s)$ w.r.t.\ the occupied parking space $x$ of the last driver by taking into account his preferred parking space $w$.\label{pic:treeparkingcases}}
\end{center}
\end{figure}

Our approach to enumerate tree parking functions for combinatorial tree families is based on a decomposition of a parking function $(T,s)$ of size $n$ w.r.t.\ the occupied parking space (i.e., parking position) $x$ of the last driver by  taking into account its preferred parking space $w=s_{n}$. Namely, when decomposing the tree $T$ into node $x$, a possibly empty sequence of subtrees $T^{(1)}, \dots, T^{(k)}$ attached to $x$ and, in case that $x$ is not the root of $T$, a subtree $T^{(0)}$, with $x$ linked to a node $y \in T^{(0)}$, one obtains that $(T^{(1)},s^{(1)}), \dots, (T^{(k)},s^{(k)})$ and, if present, $(T^{(0)},s^{(0)})$ are themselves tree parking functions, where $s^{(0)}, s^{(1)}, \dots, s^{(k)}$ are the subsequences of $s$ corresponding to the drivers arriving and parking in the respective subtree. We may distinguish the following four cases, which are illustrated in Figure~\ref{pic:treeparkingcases}:
\begin{enumerate}
\item Node $x$ is the root of $T$ and $w=x$, i.e., the last driver parks at his preferred parking space.
\item Node $x$ is the root of $T$ and $w \neq x$, i.e., the last driver does not park at his preferred parking space, thus $w \in T^{(1)} \cup T^{(2)} \cup \cdots \cup T^{(k)}$.
\item Node $x$ is not the root of $T$ and $w=x$.
\item Node $x$ is not the root of $T$ and $w \neq x$, thus $w \in T^{(1)} \cup T^{(2)} \cup \cdots \cup T^{(k)}$.
\end{enumerate}

\subsection{Ordered trees\label{ssec:parking_ordered_trees}}

Let us denote by $\mathcal{G}$ the family of tree parking functions for ordered trees. We use the decomposition of a parking function w.r.t.\ the last driver, i.e., the last occupied node in the tree, as stated in Section~\ref{ssec:ParkingFunction_Decomposition} to describe $\mathcal{G}$ via a formal recursive equation. Since the objects $g \in \mathcal{G}$ are pairs $g=(T,s)$ consisting of a tree and a sequence, we have to extend the meaning of the operator $\ast$ as a pair of operators, namely consisting of the partition product acting on the tree components (order preserving relabellings of the labels of the trees) and the shuffle product acting on the sequence components (order preserving rearrangements of the entries of the respective sequences). We define the pointing operator $\Theta$ as acting on the tree component, i.e., $\Theta(\mathcal{C})$ concerns the family of objects obtained from $\mathcal{C}$ by marking a node. We further define the boxed product as acting on the sequence components, i.e., $\mathcal{C}^{\blacksquare} \ast \mathcal{D}$ contains all objects of $\mathcal{C} \ast \mathcal{D}$, where the last driver will be assigned to the $\mathcal{C}$-component. Moreover, we may assume that the attachment points in the tree component are marked with a marker $A$; the pointing operator $\Theta_{A}$ is defined as acting on these markers in the tree component, i.e., $\Theta_{A}(\mathcal{G})$ is the family of objects obtained from $\mathcal{G}$ by distinguishing an attachment point.

Using these combinatorial operators, above decomposition immediately yields the following formal description of the family $\mathcal{G}$:
\begin{equation}
  \mathcal{G} = \mathcal{Z}^{\blacksquare} \ast \left(\Theta\big(\textsc{Seq}(\mathcal{G})\big) + \textsc{Seq}(\mathcal{G})\right) \ast \left(\{\epsilon\} + \Theta_{A}(\mathcal{G})\right),
	\label{eqn:ordered_formal_eqn_G}
\end{equation}
where $\mathcal{Z}^{\blacksquare} \ast \textsc{Seq}(\mathcal{G})$ and $\mathcal{Z}^{\blacksquare} \ast \Theta\big(\textsc{Seq}(\mathcal{G})\big)$ correspond to the situation that the last driver parks at its preferred parking space or not, respectively, and $\{\epsilon\}$ and $\Theta_{A}(\mathcal{G})$ correspond to the case that the occupied parking space of the last driver is the root of the tree or not.
According to the definition of the operator $\ast$ in the formal equation of $\mathcal{G}$ we introduce the double-exponential generating function
\begin{equation}
  G:=G(z) = \sum_{g=(T,s) \in \mathcal{G}} \frac{z^{|T|}}{(|T|!)^{2}} = \sum_{n \ge 1} G_{n} \frac{z^{n}}{(n!)^{2}},
	\label{eqn:ordered_def_G}
\end{equation}
with $G_{n} = |\{g \in \mathcal{G} : |g|=n\}|$ the total number of parking functions with $n$ drivers for size-$n$ ordered trees. We have to take into account that $\mathcal{B} = \Theta(\mathcal{C})$ corresponds to $B(z) = z C'(z)$ and that $\mathcal{B} = \mathcal{C}^{\blacksquare} \ast \mathcal{D}$ corresponds to $B'(z) = C'(z) \cdot D(z)$ at the level of generating functions (see \cite{FlaSed2009} for constructions on labelled structures, whose notion is adapted here). Moreover, we use that an ordered tree of size $n$ has exactly $2n-1$ attachment points, which shows that the generating function corresponding to the family $\Theta_{A}(\mathcal{G})$ is given by $2 z G'(z) - G(z)$. With these considerations, we obtain from \eqref{eqn:ordered_formal_eqn_G} the following first-order non-linear DEQ with initial condition $G(0)=0$:
\begin{align}
  G' & = \left(z \cdot \Big(\frac{1}{1-G}\Big)' + \frac{1}{1-G}\right) \cdot \left(1+2zG'-G\right)\notag\\
	& = \frac{1}{(1-G)^{2}} \cdot (1-G+zG') \cdot (1-G+2zG').
	\label{eqn:ordered_deq_G}
\end{align}
Below we give the computations yielding the solution of \eqref{eqn:ordered_deq_G} stated in Theorem~\ref{thm:gf_treeparking}, since first it seems that standard computer algebra systems do not directly yield the solution in this form suitable to our purpose, and second they can be adapted easily to treat tree parking functions for other combinatorial tree families.

First it seems advantageous to introduce the function $H:=H(z)$ via $H=\frac{1}{1-G}$, which implies $H' = H^{2} G'$ and yields the differential equation
\begin{equation}
  H' = (H+zH') \cdot (H+2zH').
	\label{eqn:ordered_deq_H}
\end{equation}
When setting $z=e^{t}$ and introducing $\breve{S} := \breve{S}(t) = e^{t} H(e^{t})$, which gives $H=e^{-t} \breve{S}$ and $H' = e^{-2t} (\breve{S}' - \breve{S})$, we obtain the following autonomous DEQ for $\breve{S}(t)$:
\begin{equation*}
  2 (\breve{S}')^{2} - (1+\breve{S}) \breve{S}' + \breve{S} = 0,
\end{equation*}
resp.\ (the correct branch is determined by considering $z \to 0 \Leftrightarrow t \to -\infty$, where it holds $\breve{S} \to 0$ and $\breve{S}' \to 0$):
\begin{equation}
  \breve{S}' = \frac{1+\breve{S}-\sqrt{1-6\breve{S}+\breve{S}^{2}}}{4}.
	\label{eqn:ordered_deq_S}
\end{equation}
In order to get a solution of \eqref{eqn:ordered_deq_S} in a form suitable to our purpose we consider substitutions of the following kind for $\breve{Q}:=\breve{Q}(t)$, with $a,b$ some constants to be determined:
\begin{equation*}
  \breve{Q} = \breve{S} \cdot \frac{1+a\breve{Q}}{1-b\breve{Q}},
\end{equation*}
which gives (again the correct branch is determined by considering $\breve{S} \to 0 \Leftrightarrow \breve{Q} \to 0$)
\begin{equation}
  \breve{Q} = \frac{1-a\breve{S}-\sqrt{1-(2a+4b)\breve{S}+a^{2}\breve{S}^{2}}}{2b}.
	\label{eqn:ordered_rel_S-tildeQ}
\end{equation}
Matching the square root-expressions in \eqref{eqn:ordered_deq_S} and \eqref{eqn:ordered_rel_S-tildeQ} determines $a=1$ and $b=1$, thus we carry out the substitution
\begin{equation*}
  \breve{Q} = \breve{S} \cdot \frac{1+\breve{Q}}{1-\breve{Q}}, \qquad \text{viz.} \qquad 
	\breve{Q} = \frac{1-\breve{S}-\sqrt{1-6\breve{S}+\breve{S}^{2}}}{2},
\end{equation*}
and simple computations show that
\begin{equation*}
  \breve{S}' = \frac{1-2\breve{Q}-\breve{Q}^{2}}{(1+\breve{Q})^{2}} \, \breve{Q}', \quad \text{and} \quad 
	\frac{1+\breve{S}-\sqrt{1-6\breve{S}+\breve{S}^{2}}}{4} = \frac{\breve{Q}}{1+\breve{Q}}.
\end{equation*}
Thus, when plugging the latter expressions into \eqref{eqn:ordered_deq_S}, we obtain that $\breve{Q}(t)$ satisfies the separable DEQ
\begin{equation*}
  \breve{Q}' = \frac{\breve{Q} (1+\breve{Q})}{1-2\breve{Q}-\breve{Q}^{2}},
\end{equation*}
which can be solved easily and gives, after adapting to the initial values for $z \to 0 \Leftrightarrow t \to -\infty$, the following implicit characterization of $\breve{Q}(t)$:
\begin{equation*}
  \breve{Q} = e^{t} (1+\breve{Q})^{2} e^{\breve{Q}},
\end{equation*}
respectively, when back-substituting $z=e^{t}$ and introducing the function $Q:=Q(z) = \breve{Q}(\ln z) = \breve{Q}(t)$, the characterization
\begin{equation}
  Q = z (1+Q)^{2} e^{Q}.
	\label{eqn:ordered_feq_Q}
\end{equation}
Using this defining equation \eqref{eqn:ordered_feq_Q} of $Q$, we further get the solution of \eqref{eqn:ordered_deq_H} as
\begin{equation*}
  H = \frac{\breve{S}}{z} = \frac{Q (1-Q)}{z(1+Q)} = (1-Q)(1+Q) e^{Q}.
\end{equation*}
Eventually, using $H=\frac{1}{1-G}$, we obtain the following characterization of the generating function $G=G(z)$ via the auxiliary function $Q=Q(z)$ as stated in Theorem~\ref{tab:gf_treeparking}:
\begin{equation}
  G = 1-\frac{1}{(1-Q)(1+Q)e^{Q}}, \quad \text{with} \quad Q=z(1+Q)^{2} e^{Q}.
	\label{eqn:ordered_feq_G-Q}
\end{equation}

\medskip

The explicit enumeration formul{\ae} for the different kinds of tree parking functions stated in Corollary~\ref{cor:explicitformulae_treeparking} and collected in Table~\ref{tab:explicitformulae_treeparking} can all be obtained by a straightforward application of the Lagrange-B{\"u}rmann inversion formula (see, e.g., \cite{FlaSed2009}) to the generating functions solutions stated in Theorem~\ref{thm:gf_treeparking}. Exemplarily, we give the computations for the coefficients $G_{n}$ of the g.f.\ solution \eqref{eqn:ordered_feq_G-Q} obtained right now. As it is also the case for a few other instances, the computations slightly simplify when starting with the g.f.\ $G'(z)$; simple computations show that it can be expressed by means of the auxiliary function $Q=Q(z)$ via $G'(z) = \frac{1+Q}{(1-Q)^{2}}$. With $\frac{d G'}{dQ} = \frac{3+Q}{(1-Q)^{3}}$, and $\varphi(Q) = (1+Q)^{2} e^{Q}$, we get then, for $n \ge 2$ (with $G_{1}=1$):
\begin{align*}
  [z^{n}] G(z) & = \frac{1}{n} [z^{n-1}] G'(z) = \frac{1}{n (n-1)} [Q^{n-2}] \textstyle{\frac{d G'}{dQ}} \cdot \big(\varphi(Q)\big)^{n-1} \\
	& = \frac{1}{n(n-1)} [Q^{n-2}] \frac{(3+Q) (1+Q)^{2n-2} e^{(n-1)Q}}{(1-Q)^{3}} \\
	& = \frac{1}{n(n-1)} \sum_{k=0}^{n-2} \frac{(n-1)^{k}}{k!} [Q^{n-2-k}] \frac{(3+Q) (1+Q)^{2n-2}}{(1-Q)^{3}} \\
	& = \frac{1}{n(n-1)} \sum_{k=0}^{n-2} \frac{(n-1)^{k}}{k!} \sum_{\ell=0}^{n-2-k} (\ell+1)(2\ell+3) \binom{2n-2}{n-2-k-\ell}.
\end{align*}
Using $G_{n} = (n!)^{2} \cdot [z^{n}] G(z)$, this leads to the corresponding result given in Table~\ref{tab:explicitformulae_treeparking}.

\subsection{Combinatorial tree families\label{ssec:Parking_Function_Combinatorial}}

The considerations in the previous subsection treating parking functions for ordered trees can be adapted in a rather straightforward way to other combinatorial tree families. First, the decomposition of a parking function w.r.t.\ the last driver as stated in Section~\ref{ssec:ParkingFunction_Decomposition} yields the following formal recursive description of the corresponding family $\mathcal{G}$ of parking functions:
\begin{equation}
  \mathcal{G} = \mathcal{Z}^{\blacksquare} \ast \big(\Theta(\phi(\mathcal{G})) + \phi(\mathcal{G})\big) \ast \left(\{\epsilon\} + \Theta_{A}(\mathcal{G})\right),
	\label{eqn:combtree_formal_eqn_G}
\end{equation}
where the substituted structure $\phi(\mathcal{G})$ is determined by the degree-weight generating function $\phi(t)$ stated in Table~\ref{tab:combinatorial_tree_families}. Actually, the formal equation \eqref{eqn:combtree_formal_eqn_G} holds for all simple families of labelled trees with an arbitrary $\phi(t)$, but in general the number of attachments points in a tree $T \in \mathcal{T}$ depends on the distribution of the node-degrees in $T$ and thus there is no direct link between the objects in $\mathcal{G}$ and $\Theta_{A}(\mathcal{G})$.

We introduce double-exponential generating functions
\begin{equation*}
  G:=G(z) = \sum_{g=(T,s) \in \mathcal{G}} \frac{z^{|T|}}{(|T|!)^{2}} = \sum_{n \ge 1} G_{n} \frac{z^{n}}{(n!)^{2}},
\end{equation*}
with $G_{n}$ the number of tree parking functions of size $n$, and take into account that for the combinatorial tree families considered the number of attachment points of a size-$n$ tree $T$ is a linear function of $n$ as given in Table~\ref{tab:combinatorial_tree_families}, thus at the level of generating functions $\Theta_{A}(\mathcal{G})$ corresponds to a linear combination of $G$ and $zG'$. Then \eqref{eqn:combtree_formal_eqn_G} yields the following first-order non-linear differential equations for $G(z)$:
\begin{subequations}
\begin{align}
\text{unordered trees:} & \quad  G' = \big(z e^{G} G' + e^{G}\big) (1+zG'),\\
\text{$d$-ary trees:} & \quad G' = \big(zd(1+G)^{d-1}G' + (1+G)^{d}\big) (1+G+(d-1)zG'),\\
\text{$d$-bundled trees:} & \quad G' = \Big(\frac{d z G'}{(1-G)^{d+1}} + \frac{1}{(1-G)^{d}}\Big) (1-G+(d+1)zG').
\end{align}
\label{eqn:combtree_deq_G}
\end{subequations}
Introducing the functions $H := H(z)$ via $H=\phi(G)$, above differential equations \eqref{eqn:combtree_deq_G} simplify:
\begin{subequations}
\begin{align}
\text{unordered trees:} & \quad  H' = (H+z H')^{2},\\
\text{$d$-ary trees:} & \quad H' = (H+z H') (dH + (d-1) z H'),\\
\text{$d$-bundled trees:} & \quad H' = (H+z H') (d H + (d+1) z H').
\end{align}
\label{eqn:combtree_deq_H}
\end{subequations}
These differential equations \eqref{eqn:combtree_deq_H} can be solved in a way analogous to the one for ordered trees as carried out in Section~\ref{ssec:parking_ordered_trees} and yield the results for $G(z)$ stated in Theorem~\ref{thm:gf_treeparking}. We omit these computations; however, it can be checked easily that the functions given there satisfy the differential equations \eqref{eqn:combtree_deq_G} together with the initial condition $G(0)=0$ and are thus indeed the required solutions.

\section{Parking distributions\label{sec:ParkingDistributions}}

\subsection{Combinatorial decomposition\label{ssec:ParkingDistribution_Decomposition}}

For the enumeration of tree parking distributions for combinatorial tree families one might try to use the decomposition w.r.t.\ the last driver as introduced in Section~\ref{ssec:ParkingFunction_Decomposition} for tree parking functions. However, for parking distributions $(T,\tilde{s})$, where $\tilde{s}$ might be considered as the multiset of preferred parking spaces, there is a priori no ordering on the drivers, and selecting the largest (or smallest) element of $\tilde{s}$ (i.e., to force that the last driver arrives at the parking space with largest (or smallest) label amongst the multiset of preferred parking spaces), does not seem to give tractable combinatorial descriptions. 

\begin{figure}
\begin{center}
\raisebox{1.4cm}{\includegraphics[width=1.7cm]{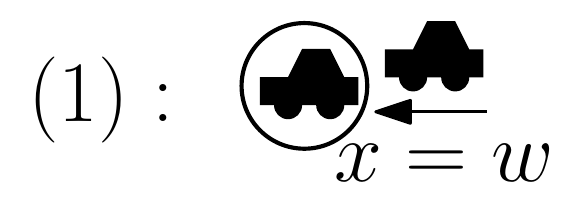}} \;
\raisebox{0.4cm}{\includegraphics[width=3.4cm]{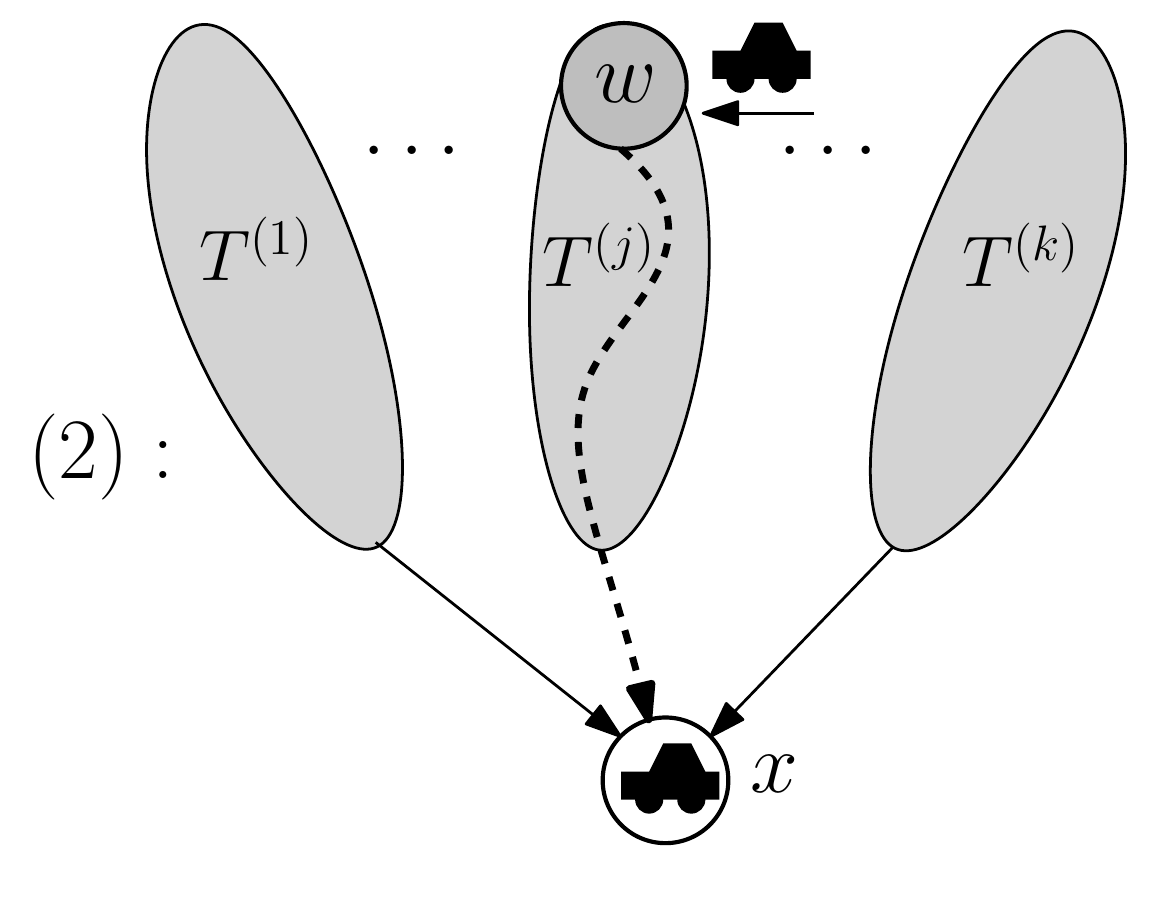}} \;
\raisebox{0.6cm}{\includegraphics[width=3.4cm]{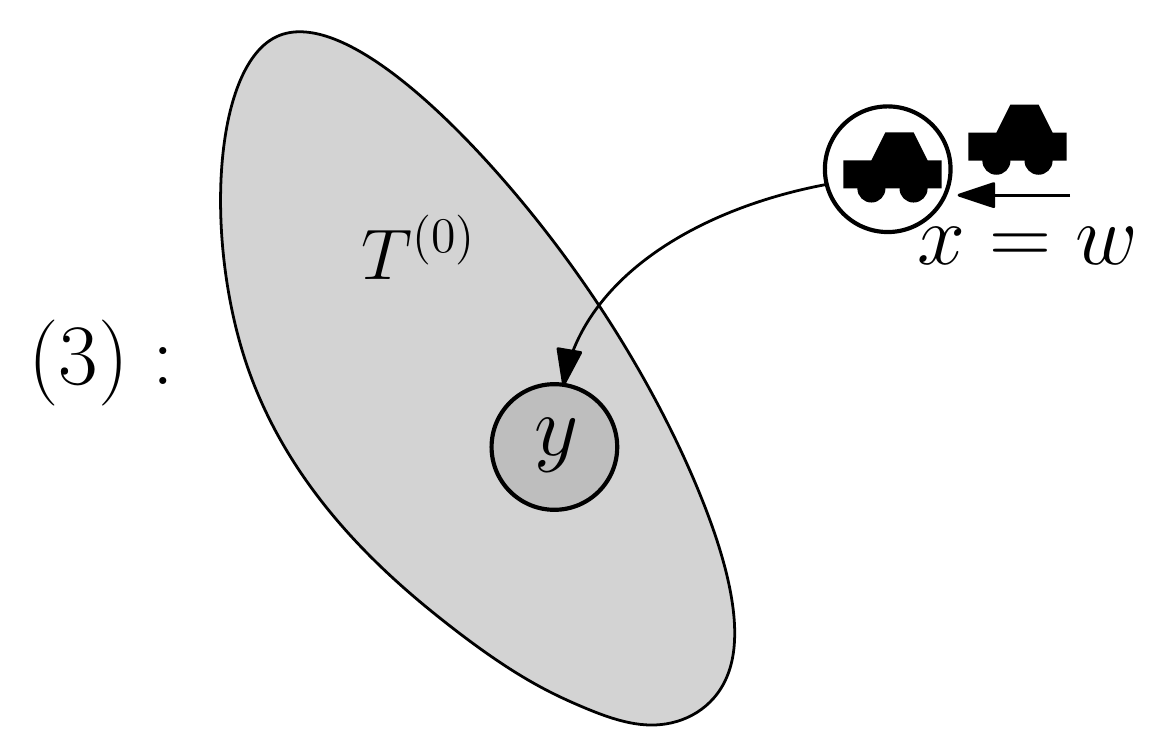}} \;
\includegraphics[width=4cm]{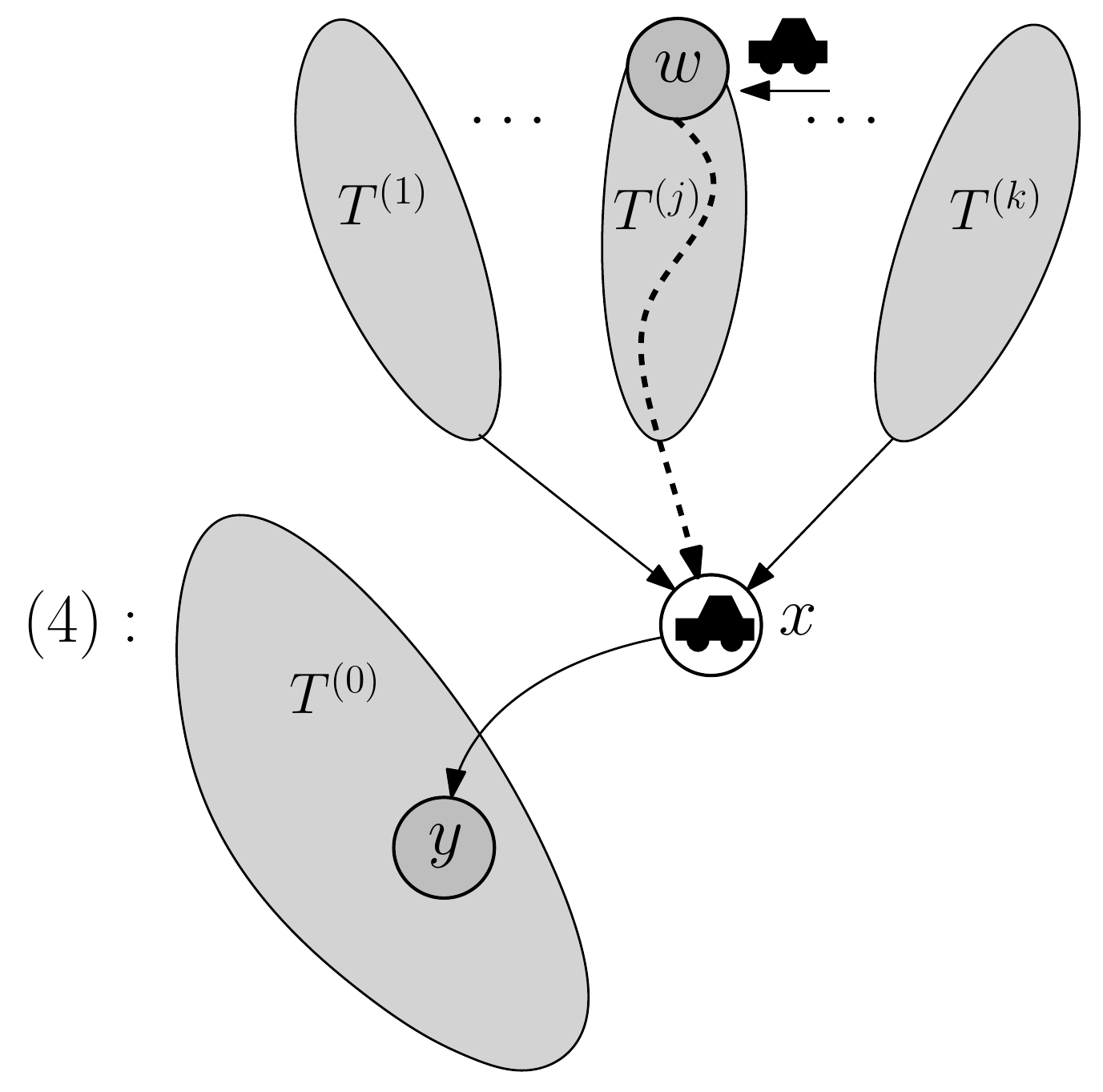}
\caption{Decomposition of a tree parking distribution $(T,\tilde{s})$ w.r.t.\ the occupied parking space $x$ of the last driver, which is assumed to arrive at a leaf $w$.\label{pic:treedistributioncases}}
\end{center}
\end{figure}

To overcome these difficulties one may assume that the last driver is arriving at a leaf of the tree; this makes sense, since for a parking function or distribution it is guaranteed that for each leaf in the tree there is at least one driver, which prefers this parking space. Thus we consider a parking distribution $(T,\tilde{s})$ and assume that the last driver arrives at leaf $w$ and parks at node $x$; then we can easily adapt the decomposition of $T$ stated in Section~\ref{ssec:ParkingFunction_Decomposition} into node $x$, a possibly empty sequence of subtrees $T^{(1)}, \dots, T^{(k)}$ attached to $x$ and, in case that $x$ is not the root of $T$, a subtree $T^{(0)}$, with $x$ linked to a node $y \in T^{(0)}$. Again, by assuming $\tilde{s}^{(0)}, \tilde{s}^{(1)}, \dots, \tilde{s}^{(k)}$ are the submultisets of $\tilde{s}$ corresponding to the drivers arriving and parking in the respective subtrees, $(T^{(1)},\tilde{s}^{(1)}), \dots, (T^{(k)},\tilde{s}^{(k)})$ and, if present, $(T^{(0)},\tilde{s}^{(0)})$ are themselves tree parking distributions. Actually, since $w$ is a leaf, some of the situations simplify and lead to the following four cases, which are illustrated in Figure~\ref{pic:treedistributioncases}:
\begin{enumerate}
\item Node $x$ is the root of $T$ and $w=x$, i.e., the last driver does not park at his preferred parking space: since $w=x$ is also a leaf, this implies $k=0$ and thus that $T$ consists of a single node.
\item Node $x$ is the root of $T$ and $w \neq x$, i.e., the last driver does not park at his preferred parking space, thus $w$ is a leaf in a tree $T^{(j)}$, $1 \le j \le k$.
\item Node $x$ is not the root of $T$ and $w=x$: again this implies $k=0$.
\item Node $x$ is not the root of $T$ and $w \neq x$, thus $w$ is a leaf in a tree $T^{(j)}$, $1 \le j \le k$.
\end{enumerate}

\subsection{Ordered trees\label{ssec:parking_distribution_ordered_trees}}

Let $\tilde{\mathcal{G}}$ be the family of tree parking distributions for ordered trees. The decomposition of a tree parking distribution stated above relies on the demand that the preferred parking space of the last driver is a leaf in the tree. To cope with this situation we have to take into account the leaves in the tree, thus we will assume that each leaf in the tree component of a tree parking distribution $(T,\tilde{s}) \in \tilde{\mathcal{G}}$ is marked with a marker $L$. Furthermore, as already used in the formal description of tree parking functions, each attachment point in the tree is marked with a marker $A$. In addition to the pointing operator $\Theta_{A}$ acting on the attachment points, we introduce the pointing operators $\Theta_{L}$ and $\tilde{\Theta}_{L}$ acting on the leaves: $\Theta_{L}(\mathcal{C})$ is the family of objects obtained from $\mathcal{C}$ by distinguishing a leaf, whereas $\tilde{\Theta}_{L}(\mathcal{C})$ is the family of objects obtained from $\mathcal{C}$ by first distinguishing a leaf and afterwards removing the respective marker $L$. Moreover, we note that, unlike in the case of parking functions, when dealing with parking distributions the operator $\ast$ is simply the partition product for labelled objects acting on the tree components, since the ordering on the entries of the multiset $\tilde{s}$ is not of relevance.

With these combinatorial operators, the decomposition described in Section~\ref{ssec:ParkingDistribution_Decomposition} leads to the following recursive description of the family $\Theta_{L}(\tilde{\mathcal{G}})$ (and thus also of $\tilde{\mathcal{G}}$):
\begin{align}
  \Theta_{L}(\tilde{\mathcal{G}}) & = \mathcal{Z} \times \{L\} + \mathcal{Z} \ast \Theta_{L}(\tilde{\mathcal{G}}) \ast \big(\textsc{Seq}(\tilde{\mathcal{G}})\big)^{2} + \mathcal{Z} \times \{L\} \ast \left(\tilde{\Theta}_{L}(\tilde{\mathcal{G}}) + \Theta_{A}(\tilde{\mathcal{G}}) - \Theta_{L}(\tilde{\mathcal{G}})\right)\notag\\
	& \quad \mbox{} + \mathcal{Z} \ast \Theta_{L}(\tilde{\mathcal{G}}) \ast \big(\textsc{Seq}(\tilde{\mathcal{G}})\big)^{2} \ast \left(\tilde{\Theta}_{L}(\tilde{\mathcal{G}}) + \Theta_{A}(\tilde{\mathcal{G}}) - \Theta_{L}(\tilde{\mathcal{G}})\right).
	\label{eqn:ordered_distribution_formal_eqn_G}
\end{align}
Here the four summands on the right hand side of \eqref{eqn:ordered_distribution_formal_eqn_G} correspond, in this order, to the four cases distinguished in the decomposition. Note that the expression $\tilde{\Theta}_{L}(\tilde{\mathcal{G}}) + \Theta_{A}(\tilde{\mathcal{G}}) - \Theta_{L}(\tilde{\mathcal{G}})$ occurring in the third and fourth summand is stemming from the distinction whether node $x$ is attached to a leaf $y \in T^{(0)}$: then the leaf in $T^{(0)}$ is no more a leaf in $T$ and thus this corresponds to $\tilde{\Theta}_{L}(\tilde{\mathcal{G}})$, or not: then a non-leaf attachment point in $T^{(0)}$ is selected and this corresponds to $\Theta_{A}(\tilde{\mathcal{G}}) - \Theta_{L}(\tilde{\mathcal{G}})$. The expression $\big(\textsc{Seq}(\tilde{\mathcal{G}})\big)^{2}$ appearing in the second and fourth summand reflects that there is a sequence of subtrees to the left and a sequence of subtrees to the right of the subtree $T^{(j)}$ containing leaf $w$.

Since this combinatorial description takes into account the leaves in the tree and according to the definition of $\ast$, we introduce bivariate (exponential, but not double-exponential) generating functions
\begin{equation}
  \check{G} := \check{G}(z,v) = \sum_{\tilde{g}=(T,\tilde{s}) \in \tilde{\mathcal{G}}} \frac{z^{|T|} \, v^{\text{$\#$ leaves of $T$}}}{|T|!}
	= \sum_{n \ge 1} \sum_{k \ge 1} \check{G}_{n,k} \frac{z^{n} v^{k}}{n!},
	\label{eqn:ordered_distribution_def_Gzv}
\end{equation}
with $\check{G}_{n,k}$ the number of tree parking distributions for ordered trees of size $n$ with $k$ leaves. Due to this definition, the families $\Theta_{L}(\tilde{\mathcal{G}})$, $\tilde{\Theta}_{L}(\tilde{\mathcal{G}})$, and $\Theta_{A}(\tilde{\mathcal{G}})$ correspond at the level of generating functions to $v \check{G}_{v}$, $\check{G}_{v}$ and $2z\check{G}_{z}-\check{G}$, respectively. Thus the formal equation \eqref{eqn:ordered_distribution_formal_eqn_G} can be translated into the following equation for $\check{G}$:
\begin{equation*}
  v \check{G}_{v} = zv + \frac{zv \check{G}_{v}}{(1-\check{G})^{2}} + zv (\check{G}_{v} + 2z\check{G}_{z}-\check{G}-v\check{G}_{v}) + \frac{zv \check{G}_{v}}{(1-\check{G})^{2}} (\check{G}_{v} + 2z\check{G}_{z} - \check{G} - v \check{G}_{v}),
\end{equation*}
from which we obtain, after simple manipulations, that $\check{G}$ satisfies the first-order non-linear PDE
\begin{equation}
  \check{G}_{v} = z \Big(1+\frac{\check{G}_{v}}{(1-\check{G})^{2}}\Big) \cdot \big(1-\check{G}+(1-v)\check{G}_{v} +2z\check{G}_{z}\big).
	\label{eqn:ordered_distribution_pde_G}
\end{equation}
Since a treatment of such PDEs seems to be less standard and we could not get solutions via standard computer algebra systems, we carry out these computations; by doing this we also state suitable constraints characterizing the particular solution of our problem.
Moreover, these considerations can be adapted easily to find also the solutions for tree parking distributions for further combinatorial tree families.

First, it turns out to be advantageous to introduce the function $H := H(z,v)$ via $H = \frac{1}{1-\check{G}}$, thus $H_{v} = H^{2} \check{G}_{v}$ and $H_{z} = H^{2} \check{G}_{z}$, which leads from \eqref{eqn:ordered_distribution_pde_G} to the following first-order non-linear PDE for $H$, which is actually the starting point of our treatment:
\begin{equation}
  H_{v} = z (1+H_{v}) \cdot (H+(1-v)H_{v}+2zH_{z}).
	\label{eqn:ordered_distribution_pde_H}
\end{equation}

In order to attack \eqref{eqn:ordered_distribution_pde_H}, we apply the method of characteristics for first-order non-linear PDEs (see, e.g., \cite{CarPea1988}). Although in general one cannot hope for getting explicit results, interestingly, for the equations occurring here, this is indeed the case.
We set $q:=H_{v}$ and $p:=H_{z}$ and rewrite \eqref{eqn:ordered_distribution_pde_H} by introducing the function $f(z,v,H,p,q)$ via
\begin{equation}
  f(z,v,H,p,q) := z (1+q) (H+(1-v)q+2zp) - q = 0.
	\label{eqn:ordered_distribution_f_Def}
\end{equation}
To get solutions of \eqref{eqn:ordered_distribution_f_Def}, we assume that each variable involved in $f$ is dependent of a parameter $t$, $z=z(t)$, $v=v(t)$, etc.{}, and study the so-called system of characteristic ordinary DEQs, which yields
\begin{equation}
\begin{split}
(i): & \quad \frac{dz}{dt} = f_{p} = 2z^{2}(1+q),\\
(ii): & \quad \frac{dv}{dt} = f_{q} = z (H+(1-v)q+2zp) + z(1+q)(1-v)-1,\\
(iii): & \quad \frac{dH}{dt} = p f_{p} + q f_{q},\\
(iv): & \quad \frac{dp}{dt} = -f_{z}-pf_{H} = -(1+q)(H+(1-v)q+2zp)-3(1+q)zp,\\
(v): & \quad \frac{dq}{dt} = -f_{v}-qf_{H} = 0.
\end{split}
\label{eqn:ordered_distribution_characteristicDEQs}
\end{equation}
We are searching for invariants, i.e., functions, where each solution curve (i.e., each characteristic curve) of the system of characteristic DEQs \eqref{eqn:ordered_distribution_characteristicDEQs} is constant. One invariant is already given via \eqref{eqn:ordered_distribution_f_Def}. Due to $(v)$, we instantly get also the invariant $q = C_{1} = \text{const}$. 
Furthermore, from \eqref{eqn:ordered_distribution_f_Def} we get $(1+q)(H + (1-v)q+2zp)=\frac{q}{z}$, and by plugging this into $(iv)$ we obtain $\frac{dp}{dt} = -\frac{q}{z} - 3(1+q)zp$. Thus, together with $(i)$ we obtain the differential equation
\begin{equation*}
  \frac{dp}{dz} = -\frac{q}{2(1+q)z^{3}} - \frac{3p}{2z}.
\end{equation*}
Taking into account that $q = C_{1}$ is constant along characteristic curves, this is a first-order linear DEQ and the general solution yields another invariant:
\begin{equation*}
  p z^{\frac{3}{2}} - \frac{q}{(1+q)\sqrt{z}} = C_{2} = \text{const}.
\end{equation*}
Moreover, again by considering \eqref{eqn:ordered_distribution_f_Def} we get $z(H+(1-v)q+2zp)=\frac{q}{1+q}$, which, by plugging into $(ii)$, leads to the equation $\frac{dv}{dt} = -\frac{1}{1+q} + z(1+q) -z(1+q)v$. Together with $(i)$ we obtain the differential equation
\begin{equation*}
  \frac{dv}{dz} = - \frac{1}{2(1+q)^{2}z^{2}} + \frac{1}{2z} - \frac{v}{2z}.
\end{equation*}
This first-order linear DEQ (again, we use that $q=C_{1}$ is constant along characteristic curves) can be solved by standard means and gives the invariant
\begin{equation*}
  \sqrt{z} (v-1) - \frac{1}{(1+q)^{2} \sqrt{z}} = C_{3} = \text{const}.
\end{equation*}
Summarizing, we get the following independent invariants of the characteristic system of DEQs \eqref{eqn:ordered_distribution_characteristicDEQs}:
\begin{equation}
\begin{split}
  q & = C_{1},\\
	pz^{\frac{3}{2}} - \frac{q}{(1+q) \sqrt{z}} & = C_{2},\\
	\sqrt{z} (v-1) - \frac{1}{(1+q)^{2} \sqrt{z}} & = C_{3},\\
	z (1+q)(H+(1-v)q+2zp)-q & = 0.
\end{split}
\label{eqn:ordered_distribution_invariants}
\end{equation}

Next we will consider the trajectories passing through $z=z_{0}$ and $v=0$, and thus evaluate the equations \eqref{eqn:ordered_distribution_invariants} at $(z,v)=(z_{0},0)$. To this aim we require an expansion of the function $H(z,v)$ around $v=0$ and thus, due to the relation $H = \frac{1}{1-\check{G}}$, an expansion of $\check{G}(z,v)$ around $v=0$, which gives the constraints required for characterizing the particular solution of our problem. According to the definition \eqref{eqn:ordered_distribution_def_Gzv}, we get $\check{G}(z,v) = v \sum_{n \ge 1} \check{G}_{n,1} \frac{z^{n}}{n!} + \mathcal{O}(v^{2})$, with $\check{G}_{n,1}$ the number of tree parking distributions of size-$n$ trees with only one leaf, which, of course, is exactly the number of parking distributions for labelled chains of length $n$. But this means that $\check{G}_{n,1} = n! \cdot B_{n}$, with $B_{n}$ the number of increasing parking functions of length $n$, for which it is well-known and already mentioned in the introduction that they are enumerated by the Catalan numbers, $B_{n} = \frac{1}{n+1} \binom{2n}{n}$, whose generating function $B:=B(z) = \sum_{n \ge 1} B_{n} z^{n}$ satisfies the functional equation $B = z (1+B)^{2}$. Thus we obtain the expansion $\check{G}(z,v) = v B(z) + \mathcal{O}(v^{2})$, and furthermore the required expansions of $H := H(z,v)$, $p := H_{z}(z,v)$ and $q := H_{v}(z,v)$ around $v=0$:
\begin{equation}
  H = \frac{1}{1-G(z,v)} = 1 + \mathcal{O}(v), \quad p=H_{z} = \mathcal{O}(v), \quad q=H_{v} = B(z) + \mathcal{O}(v).
	\label{eqn:ordered_distribution_constraints}
\end{equation}

With \eqref{eqn:ordered_distribution_constraints} we are able to determine the value of the constants $C_{j}$, $1 \le j \le 3$, of the trajectories passing through $(z_{0},0)$, for which we get
\begin{align*}
  C_{1}(z_{0}) & = q\Big|_{\substack{z=z_{0},\\ v=0}} = B(z_{0}),\\
	C_{2}(z_{0}) & = \left.pz^{\frac{3}{2}} - \frac{q}{(1+q) \sqrt{z}}\right|_{\substack{z=z_{0},\\ v=0}} = - \frac{B(z_{0})}{(1+B(z_{0})) \sqrt{z_{0}}},\\
	C_{3}(z_{0}) & = \left.\sqrt{z} (v-1) - \frac{1}{(1+q)^{2} \sqrt{z}}\right|_{\substack{z=z_{0},\\ v=0}} = -\sqrt{z_{0}} - \frac{1}{(1+B(z_{0}))^{2} \sqrt{z_{0}}}.
\end{align*}
When defining $\check{Q} := B(z_{0})$, which implies $z_{0} = \frac{B(z_{0})}{(1+B(z_{0}))^{2}} = \frac{\check{Q}}{(1+\check{Q})^{2}}$, above expressions simplify and we obtain
\begin{equation}
  C_{1}(z_{0}) = \check{Q}, \quad C_{2}(z_{0}) = -\sqrt{\check{Q}}, \quad C_{3}(z_{0}) = -\frac{1}{\sqrt{\check{Q}}}.
	\label{eqn:ordered_distribution_constants}
\end{equation}
Plugging the values \eqref{eqn:ordered_distribution_constants} into the first and third equation of \eqref{eqn:ordered_distribution_invariants}, we get $q=\check{Q}$ and further
$\sqrt{z} (v-1) - \frac{1}{(1+\check{Q})^{2} \sqrt{z}} = -\frac{1}{\sqrt{\check{Q}}}$, which, after simple manipulations, yields an implicit characterization of $\check{Q} := \check{Q}(z,v)$ via
\begin{equation}
  \check{Q} = \frac{z (1+\check{Q})^{4}}{\big(1-z (v-1) (1+\check{Q})^{2}\big)^{2}}.
	\label{eqn:ordered_distribution_Qzv}
\end{equation}
Furthermore by plugging the values \eqref{eqn:ordered_distribution_constants} into the second equation of \eqref{eqn:ordered_distribution_invariants}, we obtain from it the relation $pz^{\frac{3}{2}} = \frac{\check{Q}}{(1+\check{Q}) \sqrt{z}} - \sqrt{\check{Q}}$. Eventually, the fourth equation of \eqref{eqn:ordered_distribution_invariants} relates the function $H(z,v)$ to $\check{Q}(z,v)$ via $H=\frac{\check{Q}}{z(1+\check{Q})} -(1-v)\check{Q}-2zp = -\frac{\check{Q}}{(1+\check{Q})z} -(1-v)\check{Q} + \frac{2\sqrt{\check{Q}}}{\sqrt{z}}$, from which we get, after simple manipulations using the defining equation \eqref{eqn:ordered_distribution_Qzv} of $\check{Q}$, the following representation:
\begin{equation}
  H = \frac{(1+\check{Q})^{2} \big(1-\check{Q}-z(v-1)(1+\check{Q})^{2}\big)}{\big(1-z(v-1)(1+\check{Q})^{2}\big)^{2}}.
	\label{eqn:ordered_distribution_solution_Hzv}
\end{equation}
The function $\check{G} = \check{G}(z,v)$ is obtained from \eqref{eqn:ordered_distribution_solution_Hzv} by using $\check{G} = 1-\frac{1}{H}$, which characterizes the desired solution of \eqref{eqn:ordered_distribution_pde_G} via the auxiliary function $\check{Q}$ given in \eqref{eqn:ordered_distribution_Qzv}:
\begin{equation}
  \check{G} = 1- \frac{\big(1-z(v-1)(1+\check{Q})^{2}\big)^{2}}{(1+\check{Q})^{2} (1-\check{Q}-z(v-1)(1+\check{Q})^{2})}.
	\label{eqn:ordered_distribution_solution_Gzv}
\end{equation}

Actually, we are interested in the generating function
\begin{equation}
  \tilde{G} := \tilde{G}(z) = \sum_{\tilde{g}=(T,\tilde{s}) \in \tilde{\mathcal{G}}} \frac{z^{|T|}}{|T|!}
	= \sum_{n \ge 1} \tilde{G}_{n} \frac{z^{n}}{n!}
	\label{eqn:ordered_distribution_def_G}
\end{equation}
of the number $\tilde{G}_{n} = |\{\tilde{g} \in \tilde{\mathcal{G}} : |\tilde{g}|=n\}|$ of tree parking distributions for ordered trees of size $n$, which we simply get via $\tilde{G}(z) = \check{G}(z,1)$. Thus, by introducing the auxiliary function $Q := Q(z) = \check{Q}(z,1)$, and evaluating \eqref{eqn:ordered_distribution_Qzv} and \eqref{eqn:ordered_distribution_solution_Gzv} at $v=1$, we get the following characterization of $\tilde{G}$ stated in Theorem~\ref{thm:gf_treeparking}:
\begin{equation}
  \tilde{G} = 1-\frac{1}{(1+Q)^{2} (1-Q)}, \quad \text{with} \quad Q = z (1+Q)^{4}.
	\label{eqn:ordered_distribution_sol_G_Q}
\end{equation}

\subsection{Combinatorial tree families}

The considerations given in the previous subsection on ordered trees can be adapted to a treatment of parking distributions for further tree families. In particular, the decomposition of a tree parking distribution w.r.t.\ the last driver, which is assumed to arrive at a leaf of the tree, as given in Section~\ref{ssec:ParkingDistribution_Decomposition}, can be used to get a formal description of the family $\tilde{\mathcal{G}}$ of tree parking distributions for combinatorial tree families. Here one has to take into account also the number of attachment points at leaves in a tree $T$: for $d$-ary and $d$-bundled trees it holds that each leaf of $T$ has $d$ attachment points, whereas for unordered trees each leaf has one attachment point.

Taking into account the four cases of the combinatorial decomposition described in Section~\ref{ssec:ParkingDistribution_Decomposition} yields for combinatorial tree families $\tilde{\mathcal{G}}$ the following formal description
\begin{align}
  \Theta_{L}(\tilde{\mathcal{G}}) & = \mathcal{Z} \times \{L\} + \mathcal{Z} \ast \Theta_{L}\big(\phi(\tilde{\mathcal{G}})\big) + \mathcal{Z} \times \{L\} \ast \left(d \, \tilde{\Theta}_{L}(\tilde{\mathcal{G}}) + \Theta_{A}(\tilde{\mathcal{G}}) - d \, \Theta_{L}(\tilde{\mathcal{G}})\right)\notag\\
	& \quad \mbox{} + \mathcal{Z} \ast \Theta_{L}\big(\phi(\tilde{\mathcal{G}})\big) \ast \left(d \, \tilde{\Theta}_{L}(\tilde{\mathcal{G}}) + \Theta_{A}(\tilde{\mathcal{G}}) - d \, \Theta_{L}(\tilde{\mathcal{G}})\right),
	\label{eqn:combinatorialtree_distribution_formal_eqn_G}
\end{align}
where we set $d=1$ for unordered trees and $\phi$ is given by the degree-weight generating function $\phi(t)$ stated in Table~\ref{tab:combinatorial_tree_families}.

Introducing bivariate generating functions
\begin{equation}
  \check{G} := \check{G}(z,v) = \sum_{\tilde{g}=(T,\tilde{s}) \in \tilde{\mathcal{G}}} \frac{z^{|T|} \, v^{\text{$\#$ leaves of $T$}}}{|T|!}
	= \sum_{n \ge 1} \sum_{k \ge 1} \check{G}_{n,k} \frac{z^{n} v^{k}}{n!},
	\label{eqn:combinatorialtree_distribution_def_Gzv}
\end{equation}
with $\check{G}_{n,k}$ the number of tree parking distributions for trees of size $n$ with $k$ leaves gives, by taking into account the number of attachment points for the corresponding tree families, the following first-order non-linear PDEs for $\check{G}(z,v)$:
\begin{subequations}
\begin{align}
\text{unordered trees:} & \quad  \check{G}_{v} = z \big(1+e^{\check{G}} \check{G}_{v}\big) (1 + (1-v)\check{G}_{v} + z\check{G}_{z}),\\
\text{$d$-ary trees:} & \quad \check{G}_{v} = z \big(1+d(1+\check{G})^{d-1}\check{G}_{v}\big) (1+\check{G} + d(1-v)\check{G}_{v} + (d-1)z\check{G}_{z}),\\
\text{$d$-bundled trees:} & \quad \check{G}_{v} = z \Big(1+\frac{d}{(1-\check{G})^{d+1}} \check{G}_{v}\Big) (1-\check{G}+d(1-v)\check{G}_{v} + (d+1)z\check{G}_{z}).
\end{align}
\label{eqn:combinatorialtree_distribution_pde_G}
\end{subequations}
Introducing the functions $H := H(z,v)$ via $H=\phi(\check{G})$, above PDEs \eqref{eqn:combinatorialtree_distribution_pde_G} simplify:
\begin{subequations}
\begin{align}
\text{unordered trees:} & \quad  H_{v} = z (1+H_{v})(H + (1-v)H_{v} + zH_{z}),\\
\text{$d$-ary trees:} & \quad H_{v} = z (1+H_{v}) (d H + d(1-v)H_{v} + (d-1)zH_{z}),\\
\text{$d$-bundled trees:} & \quad H_{v} = z (1+H_{v}) (d H + d(1-v)H_{v} + (d+1)zH_{z}).
\end{align}
\label{eqn:combinatorialtree_distribution_pde_H}
\end{subequations}
These differential equations \eqref{eqn:combinatorialtree_distribution_pde_H} can be solved by the method of characteristics analogous to the computations for ordered trees carried out in Section~\ref{ssec:parking_distribution_ordered_trees} and lead to the following solutions for $\check{G} := \check{G}(z,v)$ dependent on respective auxiliary functions $\check{Q} := \check{Q}(z,v)$:
\begin{subequations}
\begin{align}
\begin{split}
\text{unordered trees:} & \quad \check{G} = \ln\left((1-\check{Q})(1+\check{Q}) e^{\check{Q}+(v-1)z(1+\check{Q})^{2}} + (1-v)\check{Q}^{2}\right),\\
& \qquad \qquad \text{with} \quad \check{Q} = z(1+\check{Q})^{2} e^{\check{Q}+(v-1)z(1+\check{Q})^{2}},
\end{split}\\
\begin{split}
\text{$d$-ary trees:} & \quad {\textstyle{\check{G} = \frac{(1+\check{Q})^{\frac{1}{d}} \big(1-\check{Q}+(v-1)z(1+\check{Q})(1-d\check{Q})\big)^{\frac{1}{d}} \big(1+(v-1)z(1+\check{Q})^{2}\big)^{\frac{d-1}{d}}}{1-\frac{\check{Q}}{d}} - 1}},\\
& \qquad \qquad \text{with} \quad \check{Q} = \frac{d z (1+\check{Q})^{2} \big(1+(v-1)z(1+\check{Q})^{2}\big)^{d-1}}{(1-\frac{\check{Q}}{d})^{d-1}},
\end{split}\\
\begin{split}
\text{$d$-bundled trees:} & \quad \check{G} = 1 - \frac{\big(1-(v-1)z(1+\check{Q})^{2}\big)^{\frac{d+1}{d}}}{(1+\frac{\check{Q}}{d})(1+\check{Q})^{\frac{1}{d}} \big(1-\check{Q}-(v-1)z(1+\check{Q})(1+d\check{Q})\big)^{\frac{1}{d}}},\\
& \qquad \qquad \text{with} \quad \check{Q} = \frac{d z (1+\check{Q})^{2} (1+\frac{\check{Q}}{d})^{d+1}}{\big(1-(v-1)z(1+\check{Q})^{2}\big)^{d+1}}.
\end{split}
\end{align}
\label{eqn:combinatorialtree_distribution_sol_Gzv}
\end{subequations}
We omit the computations, but remark that it can be checked rather easily that these functions indeed satisfy the PDEs \eqref{eqn:combinatorialtree_distribution_pde_G}.

Introducing the generating functions
\begin{equation}
  \tilde{G} := \tilde{G}(z) = \sum_{\tilde{g}=(T,\tilde{s}) \in \tilde{\mathcal{G}}} \frac{z^{|T|}}{|T|!}
	= \sum_{n \ge 1} \tilde{G}_{n} \frac{z^{n}}{n!}
	\label{eqn:combinatorialtree_distribution_def_G}
\end{equation}
of the number $\tilde{G}_{n} = |\{\tilde{g} \in \tilde{\mathcal{G}} : |\tilde{g}|=n\}|$ of tree parking distributions of size $n$ for the respective combinatorial tree family, we obtain via $\tilde{G}(z) = \check{G}(z,1)$ and by introducing the auxiliary functions $Q := Q(z) = \check{Q}(z,1)$ from \eqref{eqn:combinatorialtree_distribution_sol_Gzv} the characterizations of $\tilde{G}$ stated in Theorem~\ref{thm:gf_treeparking} and collected in Table~\ref{tab:gf_treeparking}.

\subsection{Extensions to further tree families\label{ssec:parking_distribution_further_trees}}

Although in this paper we mainly restrict ourselves to combinatorial tree families as introduced in Section~\ref{ssec:Tree_families}, the approach is not limited to these tree families, but could be extended further. In particular, for tree parking distributions, where also the leaves are taken into account in the combinatorial decomposition of the trees, the considerations of this section can be adapted without further difficulties to such simple tree families, where the number of attachment points $A(T)$ in a tree $T$ depends only on the size and the number of leaves of the tree. Amongst them are various tree families occurring in the combinatorial literature as ordered unary-binary trees (also called Motzkin-trees) with degree-weight g.f.\ $\phi(t) = 1+t+t^{2}$, unordered unary-binary trees with degree-weight g.f.\ $\phi(t) = 1 + t + \frac{t^{2}}{2}$, and mobile trees with degree-weight g.f.\ $\phi(t) = 1 + \ln(\frac{1}{1-t})$. Another family that can be treated rather easily are strict binary trees with degree-weight g.f.\ $\phi(t) = 1+t^{2}$, since we can use the solution for the bivariate g.f.\ $\check{G}(z,v)$ of binary trees obtained from \eqref{eqn:combinatorialtree_distribution_sol_Gzv} and consider $\left.v \check{G}(zv,v^{-2})\right|_{v=0}$, which gives exactly the g.f.\ of parking distributions for strict binary trees, i.e., binary trees without nodes of in-degree $1$. We state results for the generating functions $\tilde{G}:=\tilde{G}(z) = \sum_{n \ge 1} \tilde{G}_{n} \frac{z^{n}}{n!}$ of the number $\tilde{G}_{n}$ of tree parking distributions of size $n$ for the aforementioned tree families in Table~\ref{tab:gf_further_treeparking}, but only exemplify a sketch of the computations for Motzkin trees.

\begin{table}
\renewcommand{\arraystretch}{1.5}
\begin{equation*}
\begin{array}{c|c|c}
\text{tree family} & \text{g.f.\ $\tilde{G}=\tilde{G}(z)$} & \text{auxiliary g.f.\ $Q=Q(z)$}\\
\hline \hline
\rule[-2ex]{0mm}{5ex}\text{Motzkin trees} & \tilde{G} = \frac{(1+2Q)\sqrt{\frac{1+Q-3Q^{2}}{1+Q+Q^{2}}}-1}{2} & Q = \frac{z \big(1+2Q+3Q^{2}\big)^{2}}{1+Q+Q^{2}}\\
\hline
\text{unordered unary-binary trees} & \tilde{G} = (1+Q) \sqrt{\frac{1+Q-\frac{3Q^{2}}{2}}{1+Q+\frac{Q^{2}}{2}}}-1 & Q = \frac{z \big(1+2Q+\frac{3Q^{2}}{2}\big)^{2}}{1+Q+\frac{Q^{2}}{2}}\\
\hline
\text{mobile trees} & \tilde{G} = 1-\frac{1}{1+Q} \, e^{\frac{Q^{2}}{1-\ln(\frac{1}{1+Q})}} & Q = \frac{z(1+Q) \big(1+Q-\ln(\frac{1}{1+Q})\big)^{2}}{1-\ln(\frac{1}{1+Q})}\\
\hline
\text{strict binary trees} & \tilde{G} = \sqrt{\frac{Q(1-2Q)}{2-Q}} & Q = \frac{4z^{2} (1+Q)^{4}}{2-Q}\\
\hline
\end{array}
\end{equation*}
\caption{Generating functions solutions $\tilde{G}$ of tree parking distributions for further tree families.\label{tab:gf_further_treeparking}}
\end{table}

Considering a Motzkin-tree $T$ with $n_{0}$, $n_{1}$, and $n_{2}$ nodes of in-degree $0$, $1$, and $2$, respectively, the number of attachment points $A(T)$ is given by $n_{0}+2n_{1}$. Thus $A(T)=2n+2-3k$ for a tree $T$ of size $n:=|T|=n_{0}+n_{1}+n_{2}$ with $k:=n_{0}$ leaves, with $k$ attachment points at leaves and $2n+2-4k$ attachment points at non-leaves.
Thus, when taking into account the four cases of the combinatorial decomposition described in Section~\ref{ssec:ParkingDistribution_Decomposition}, we obtain for the family $\tilde{\mathcal{G}}$ of tree parking distributions for Motzkin-trees the formal description (with $\phi(t)=1+t+t^{2}$):
\begin{align}
  \Theta_{L}(\tilde{\mathcal{G}}) & = \mathcal{Z} \times \{L\} + \mathcal{Z} \ast \Theta_{L}\big(\phi(\tilde{\mathcal{G}})\big) + \mathcal{Z} \times \{L\} \ast \left(\tilde{\Theta}_{L}(\tilde{\mathcal{G}}) + 2 \Theta(\tilde{\mathcal{G}}) + 2 \tilde{\mathcal{G}} -4 \Theta_{L}(\tilde{\mathcal{G}})\right)\notag\\
	& \quad \mbox{} + \mathcal{Z} \ast \Theta_{L}\big(\phi(\tilde{\mathcal{G}})\big) \ast \left(\tilde{\Theta}_{L}(\tilde{\mathcal{G}}) + 2 \Theta(\tilde{\mathcal{G}}) + 2\tilde{\mathcal{G}} -4\Theta_{L}(\tilde{\mathcal{G}})\right).
	\label{eqn:Motzkin_distribution_formal_eqn_G}
\end{align}

Introducing bivariate generating functions
\begin{equation}
  \check{G} := \check{G}(z,v) = \sum_{n \ge 1} \sum_{k \ge 1} \check{G}_{n,k} \frac{z^{n} v^{k}}{n!},
	\label{eqn:Motzkin_distribution_def_Gzv}
\end{equation}
with $\check{G}_{n,k}$ the number of tree parking distributions of trees of size $n$ with $k$ leaves gives, after simple manipulations, the following first-order non-linear PDE for $\check{G}(z,v)$:
\begin{equation*}
  \check{G}_{v} = z \big(1+(1+2\check{G}) \check{G}_{v}\big) (1 + 2\check{G} + (1-4v)\check{G}_{v} + 2z\check{G}_{z}).
\end{equation*}
Introducing the function $\check{H} = \check{H}(z,v)$ via $\check{H}=1+\check{G}+\check{G}^{2}$ the PDE slightly simplifies and yields
\begin{equation*}
  \check{H}_{v} = z (1+\check{H}_{v}) (4\check{H} -3 +(1-4v) \check{H}_{v} +2z\check{H}_{z}).
\end{equation*}
An application of the method of characteristics carried out as in Section~\ref{ssec:parking_distribution_ordered_trees} yields the following solution for $\check{H} := \check{H}(z,v)$ (and thus also for $\check{G}$) dependent on an auxiliary function $\check{N} := \check{N}(z,v)$:
\begin{equation*}
  \check{H} = \frac{\check{N}}{2z(1+\check{N})} - \frac{\check{N}^{2}}{4z^{2}(1+\check{N})^{3}} + \frac{3+(4v-1)\check{N}}{4}, \qquad \check{N} = \frac{z(1+\check{N})^{2} (2+z(4v-1)(1+\check{N})^{2})}{2-\check{N}}.
\end{equation*}

Since we are only interested in the generating function
\begin{equation*}
  \tilde{G} := \tilde{G}(z) = \sum_{n \ge 1} \tilde{G}_{n} \frac{z^{n}}{n!}
\end{equation*}
of the number $\tilde{G}_{n}$ of tree parking distributions of size $n$, which we simply get via $\tilde{G}(z) = \check{G}(z,1)$, we also introduce the auxiliary functions $H := H(z) = \check{H}(z,1)$ and $N := N(z) = \check{N}(z,1)$ and evaluate above representation at $v=1$. This yields the following implicit characterization of $\tilde{G}$:
\begin{equation*}
  H = 1 + \tilde{G} + \tilde{G}^{2}, \quad H=\frac{(1-N)(N-z(1+N)^{2})}{2z^{2}(1+N)^{3}}, \quad N = \frac{2z(1+N)^{2}}{2-N} + \frac{3z^{2}(1+N)^{4}}{2-N}.
\end{equation*}
However, it is appropriate to introduce another auxiliary function $Q := Q(z)$ via the substitution $N = \frac{Q(1+2Q)}{1+Q+Q^{2}}$, which eventually leads to the result stated in Table~\ref{tab:gf_further_treeparking}.

\section{Prime tree parking functions and distributions\label{sec:Prime}}

\subsection{Combinatorial decomposition\label{ssec:PrimeParking_Decomposition}}

\begin{figure}
\begin{center}
\includegraphics[width=4.5cm]{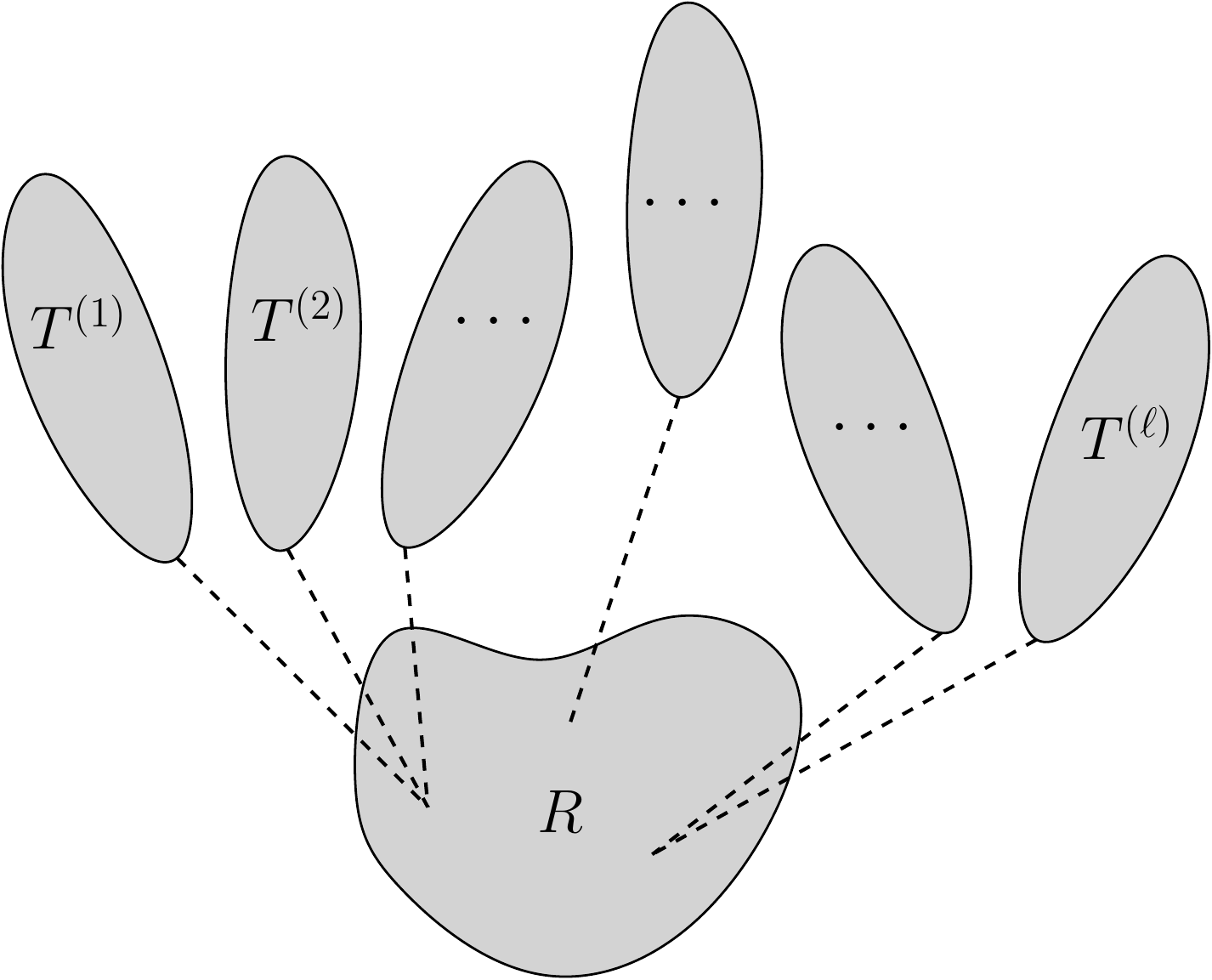}
\caption{Decomposition of a tree parking function $(T,s)$ w.r.t.\ the root-core $R$. The subtrees $T^{(j)}$, $1 \le j \le \ell$, are attached to nodes of $R$ via unused edges (visualized by dashed edges).\label{pic:primetreedecomposition}}
\end{center}
\end{figure}

The enumeration of prime tree parking functions and distributions can be obtained from the results gained in Section~\ref{sec:ParkingFunctions}-\ref{sec:ParkingDistributions} by using the decomposition of a (arbitrary) parking function w.r.t.\ the so-called root-core of the tree component. Namely, for a tree parking function $(T,s)$, let us denote by $R$ the maximal subtree of $T$ containing the root of $T$ and only used edged (i.e., edges that are used by drivers during the parking procedure). Removing $R$, $T$ decomposes into trees $T^{(1)}, \dots, T^{(\ell)}$, where each of these subtrees is originally attached to a node of $R$ via an unused edge (i.e., one that is not used by any of the drivers during the parking procedure). When we denote by $p$ and $s^{(1)}, \dots, s^{(\ell)}$ the subsequences of $s$ corresponding to the drivers arriving (and parking) in the subtrees $R$ and $T^{(1)}, \dots, T^{(\ell)}$, respectively, then it holds that $(R,p)$ is a prime tree parking function and $(T^{(1)},s^{(1)}), \dots, (T^{(\ell)},s^{(\ell)})$ are (arbitrary) tree parking functions. The decomposition is illustrated in Figure~\ref{pic:primetreedecomposition}. Of course, the decomposition can be used in a completely analogous way also for tree parking distributions.

\subsection{Ordered trees}

\subsubsection{Prime tree parking functions}
Here we study prime tree parking functions and distributions for ordered trees and start with parking functions. Thus let us denote by $\mathcal{P}$ the family of prime tree parking functions and by $\mathcal{G}$ the family of (arbitrary) tree parking functions for ordered trees. We may assume that for any prime tree parking function each attachment point in the tree component is marked with a marker $A$. We may express via $\mathcal{P} = \mathcal{P}(\mathcal{Z},A)$ the possession of nodes $\mathcal{Z}$ and attachment points $A$ for each object of $\mathcal{P}$. The decomposition of a tree parking function w.r.t.\ the root-core as given in the previous subsection yields the formal equation
\begin{equation}
  \mathcal{G} = \mathcal{P}\big(\mathcal{Z}, \textsc{Seq}(\mathcal{G})\big),
	\label{eqn:prime_parking_formal_P-G}
\end{equation}
since at each attachment point in the root-core (which actually corresponds to a prime tree parking function) there is attached a (possibly empty) sequence of (arbitrary) tree parking functions, thus can be described via a substitution of $A$ by $\textsc{Seq}(\mathcal{G})$.

Introducing the (double-exponential) generating functions
\begin{equation*}
  P:=P(z) = \sum_{p=(T,s) \in \mathcal{P}} \frac{z^{|T|}}{(|T|!)^{2}} = \sum_{n \ge 1} P_{n} \frac{z^{n}}{(n!)^{2}},
\end{equation*}
with $P_{n} = |\{p \in \mathcal{P} : |p|=n\}|$ the total number of prime tree parking functions of size $n$ for ordered trees, using the generating function $G(z)$ of the number of (arbitrary) tree parking functions as defined in \eqref{eqn:ordered_def_G}, and taking into account that any ordered tree $T$ of size $n$ has exactly $2n-1$ attachment points, the formal equation~\eqref{eqn:prime_parking_formal_P-G} leads at the level of generating functions to the following relation between $G = G(z)$ and $P = P(z)$:
\begin{equation}
  G = (1-G) \cdot P\left(\frac{z}{(1-G)^{2}}\right), \qquad \text{i.e.}, \qquad  \frac{G}{1-G} = P\left(\frac{z}{(1-G)^{2}}\right).
	\label{eqn:prime_parking_eqn_P-G}
\end{equation}

Now we use the characterization \eqref{eqn:ordered_feq_G-Q} of the generating function $G$ by means of an auxiliary function $Q$, and express the left hand-side and the right hand-side, respectively, of equation \eqref{eqn:prime_parking_eqn_P-G} via this function $Q$ by using its defining functional equation. We obtain
\begin{equation*}
  \frac{G}{1-G} = (1+Q)(1-Q)e^{Q}-1 \quad \text{and} \quad \frac{z}{(1-G)^{2}} = Q (1-Q)^{2} e^{Q},
\end{equation*}
thus the function $P$ is determined via
\begin{equation}
  P\big(Q (1-Q)^{2} e^{Q}\big) = (1+Q)(1-Q)e^{Q}-1.
	\label{eqn:prime_parking_sol_P-Q}
\end{equation}
Therefore we introduce another auxiliary function $Q := Q(z)$ that is characterized implicitly via the functional equation
\begin{equation}
  z = Q (1-Q)^{2} e^{Q}, \qquad \text{resp.} \qquad Q = \frac{z}{(1-Q)^{2} \, e^{Q}},
	\label{eqn:prime_parking_sol_Q}
\end{equation}
such that the desired generating function $P = P(z)$ is obtained from it via
\begin{equation}
  P = (1+Q)(1-Q)e^{Q}-1,
	\label{eqn:prime_parking_sol_P}
\end{equation}
which is exactly the result stated in Theorem~\ref{thm:gf_treeparking}.

\subsubsection{Prime tree parking distributions}
All these considerations can be applied in an analogous way also for prime tree parking distributions. Let $\tilde{\mathcal{P}}$ and $\tilde{\mathcal{G}}$ denote the family of prime parking distributions and (arbitrary) tree parking distributions, respectively, for ordered trees. Again, by expressing via $\tilde{\mathcal{P}} = \tilde{\mathcal{P}}(\mathcal{Z},A)$ that any object of $\tilde{\mathcal{P}}$ possesses nodes $\mathcal{Z}$ and attachment points $A$, the decomposition of a tree parking distribution w.r.t.\ the root-core leads to an analogue of the formal equation \eqref{eqn:prime_parking_formal_P-G} relating the families $\tilde{\mathcal{P}}$ and $\tilde{\mathcal{G}}$.

Introducing the exponential generating function
\begin{equation*}
  \tilde{P} := \tilde{P}(z) = \sum_{\tilde{p}=(T,\tilde{s}) \in \tilde{\mathcal{P}}} \frac{z^{|T|}}{|T|!} = \sum_{n \ge 1} \tilde{P}_{n} \frac{z^{n}}{n!},
\end{equation*}
with $\tilde{P}_{n} = |\{\tilde{p} \in \tilde{\mathcal{P}} : |\tilde{p}|=n\}|$ the total number of prime tree parking distributions of size $n$ for ordered trees, and using the generating function $\tilde{G}(z)$ of the number of (arbitrary) tree parking distributions as defined in \eqref{eqn:ordered_distribution_def_G}, we obtain from the formal equation~\eqref{eqn:prime_parking_formal_P-G} the following relation between $\tilde{G} = \tilde{G}(z)$ and $\tilde{P} = \tilde{P}(z)$:
\begin{equation}
  \frac{\tilde{G}}{1-\tilde{G}} = \tilde{P}\left(\frac{z}{(1-\tilde{G})^{2}}\right).
	\label{eqn:prime_distribution_eqn_P-G}
\end{equation}

Using the characterization \eqref{eqn:ordered_distribution_sol_G_Q} of the generating function $\tilde{G}$ by means of an auxiliary function $Q$ and expressing the left hand-side and the right hand-side, respectively, of equation \eqref{eqn:prime_distribution_eqn_P-G} via this function $Q$ by using its defining functional equation, we obtain for $\tilde{P}$:
\begin{equation}
  \tilde{P}\big(Q (1-Q)^{2}\big) = (1+Q)^{2}(1-Q)-1.
	\label{eqn:prime_distribution_sol_P-Q}
\end{equation}
This yields the characterization of $\tilde{P} = \tilde{P}(z)$ via an auxiliary function $Q = \tilde{Q}(z)$ as stated in Theorem~\ref{thm:gf_treeparking}:
\begin{equation}
  \tilde{P} = (1+Q)^{2}(1-Q)-1, \qquad \text{with} \qquad Q = \frac{z}{(1-Q)^{2}}.
	\label{eqn:prime_distribution_sol_P}
\end{equation}

\subsection{Combinatorial tree families}

The considerations yielding enumerative results for prime tree parking functions and distributions from the corresponding findings for (arbitrary) tree parking functions and distributions, respectively, just carried out for ordered trees can be extended to the other combinatorial tree families. In the following we only sketch them for parking functions.

The decomposition of a tree parking function w.r.t.\ the root-core can be translated into a formal equation relating the family of tree parking functions $\mathcal{G}$ and the family of prime tree parking functions $\mathcal{P} = \mathcal{P}(\mathcal{Z},A)$, where one has to take into account that for unordered trees at each attachment point in the root-core there is attached a (possibly empty) set of (arbitrary) tree parking functions, whereas for $d$-ary trees at each attachment point of the root-core there is either attached a tree parking function or not. Furthermore, as for the particular instance of ordered trees, for $d$-bundled trees at each attachment point in the root-core there is attached a sequence of parking functions. This gives the following symbolic equations:
\begin{align}
\text{unordered trees:} & \quad \mathcal{G} = \mathcal{P}\big(\mathcal{Z},\textsc{Set}(\mathcal{G})\big),\notag\\
\text{$d$-ary trees:} & \quad \mathcal{G} = \mathcal{P}\big(\mathcal{Z},\{\epsilon\}+\mathcal{G}\big),\label{eqn:general_prime_parking_formal_P-G}\\
\text{$d$-bundled trees:} & \quad \mathcal{G} = \mathcal{P}\big(\mathcal{Z},\textsc{Seq}(\mathcal{G})\big).\notag
\end{align}

Introducing the generating functions
\begin{equation*}
  P:=P(z) = \sum_{p=(T,s) \in \mathcal{P}} \frac{z^{|T|}}{(|T|!)^{2}} = \sum_{n \ge 1} P_{n} \frac{z^{n}}{(n!)^{2}},
\end{equation*}
with $P_{n} = |\{p \in \mathcal{P} : |p|=n\}|$ the total number of prime tree parking functions of size $n$ for the corresponding tree family, and keeping in mind the number of attachment points of a size-$n$ tree in the respective tree families stated in Table~\ref{tab:combinatorial_tree_families}, we get from \eqref{eqn:general_prime_parking_formal_P-G} the following relations between the generating functions $G := G(z)$ (of the number of arbitrary tree parking functions defined in Section~\ref{ssec:Parking_Function_Combinatorial}) and $P := P(z)$:
\begin{align}
\text{unordered trees:} & \quad G = P\left(z e^{G}\right),\notag\\
\text{$d$-ary trees:} & \quad G = (1+G) \cdot P\left(z (1+G)^{d-1}\right),\label{eqn:general_prime_parking_eqn_P-G}\\
\text{$d$-bundled trees:} & \quad G = (1-G) \cdot P\left(\frac{z}{(1-G)^{d+1}}\right).\notag
\end{align}
Computations analogous to the ones for ordered trees, by taking into account the characterization of the g.f.\ for $G(z)$ obtained previously and collected in Table~\ref{tab:gf_treeparking}, show the results stated in Theorem~\ref{thm:gf_treeparking}. Again it is not difficult to show that the functions given there satisfy above equations \eqref{eqn:general_prime_parking_eqn_P-G} and are thus indeed the required solutions.

\medskip

We remark that for the particular instance of prime tree parking distributions for binary trees an application of the Lagrange-B{\"u}rmann inversion formula to the g.f.\ solution would not directly yield the explicit formula stated in Table~\ref{tab:explicitformulae_treeparking}. Instead, one can easily get it from the representation 
$\tilde{P} = 1-\frac{2-Q}{2 \, \sqrt{1-Q^{2}}}$, $Q = \frac{2z (1+Q)^{\frac{3}{2}}}{\sqrt{1-Q}}$ given in Table~\ref{tab:gf_treeparking} as follows. Simple computations show that $\tilde{P}'(z) = \frac{1+Q}{1-Q}$, and from the functional equation of $Q$ we get
\begin{equation*}
  \frac{1+Q}{1-Q} = \frac{1}{1-Q} + 2z \left(\frac{1+Q}{1-Q}\right)^{\frac{3}{2}}.
\end{equation*}
Since $\frac{1}{1-Q} = \frac{1}{2} \left(1+\frac{1+Q}{1-Q}\right)$, we obtain $\tilde{P}' = \frac{1}{2} (1+\tilde{P}') + 2z (\tilde{P}')^{\frac{3}{2}}$, resp.\
\begin{equation*}
  \tilde{P}' = 1 + 4z \big(\tilde{P}'\big)^{\frac{3}{2}}.
\end{equation*}
Setting $A=\tilde{P}'-1$ and applying the Lagrange-B{\"u}rmann inversion formula to the resulting equation $A = 4z (1+A)^{\frac{3}{2}}$ shows the result for $\tilde{P}_{n}$ stated in Table~\ref{tab:explicitformulae_treeparking}.

\section{General tree parking functions and distributions\label{sec:General}}

\subsection{Combinatorial decomposition\label{ssec:GeneralParking_Decomposition}}

\begin{figure}
\begin{center}
\raisebox{5mm}{\includegraphics[width=4cm]{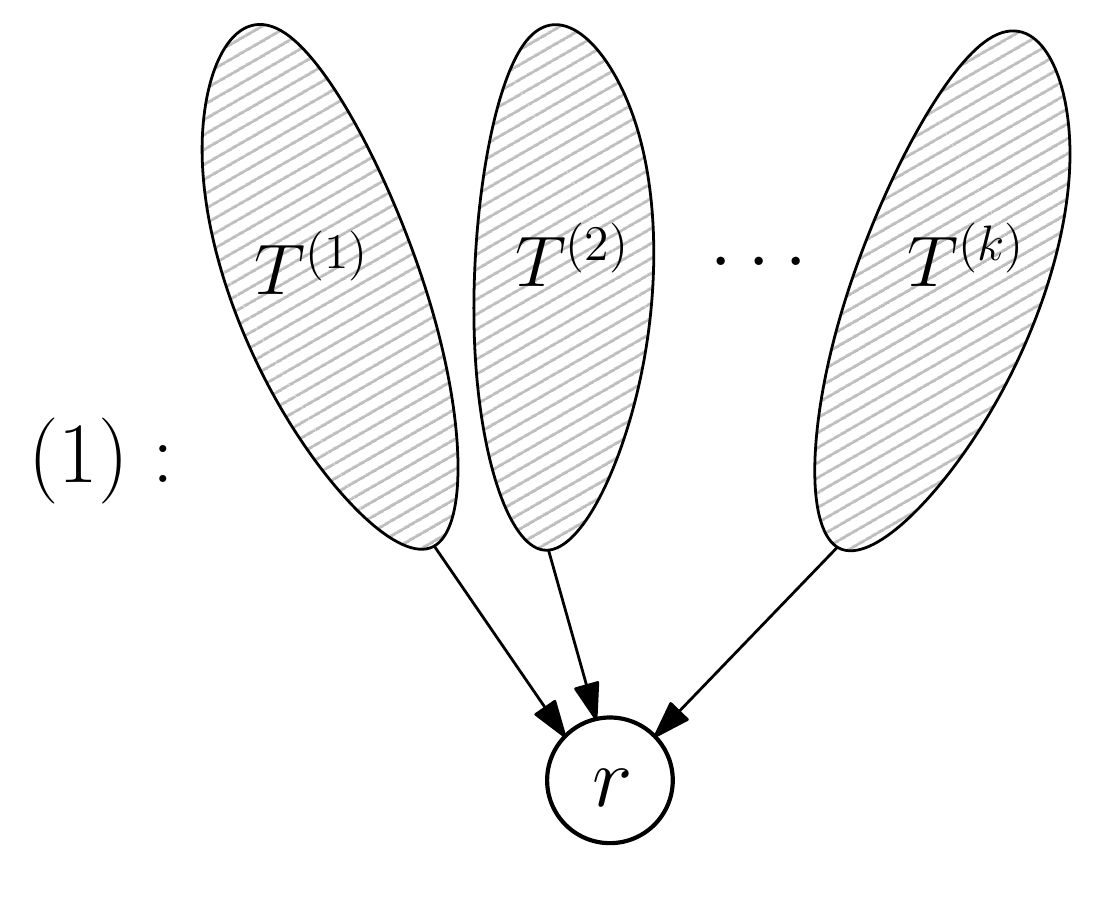}} \qquad
\includegraphics[width=7cm]{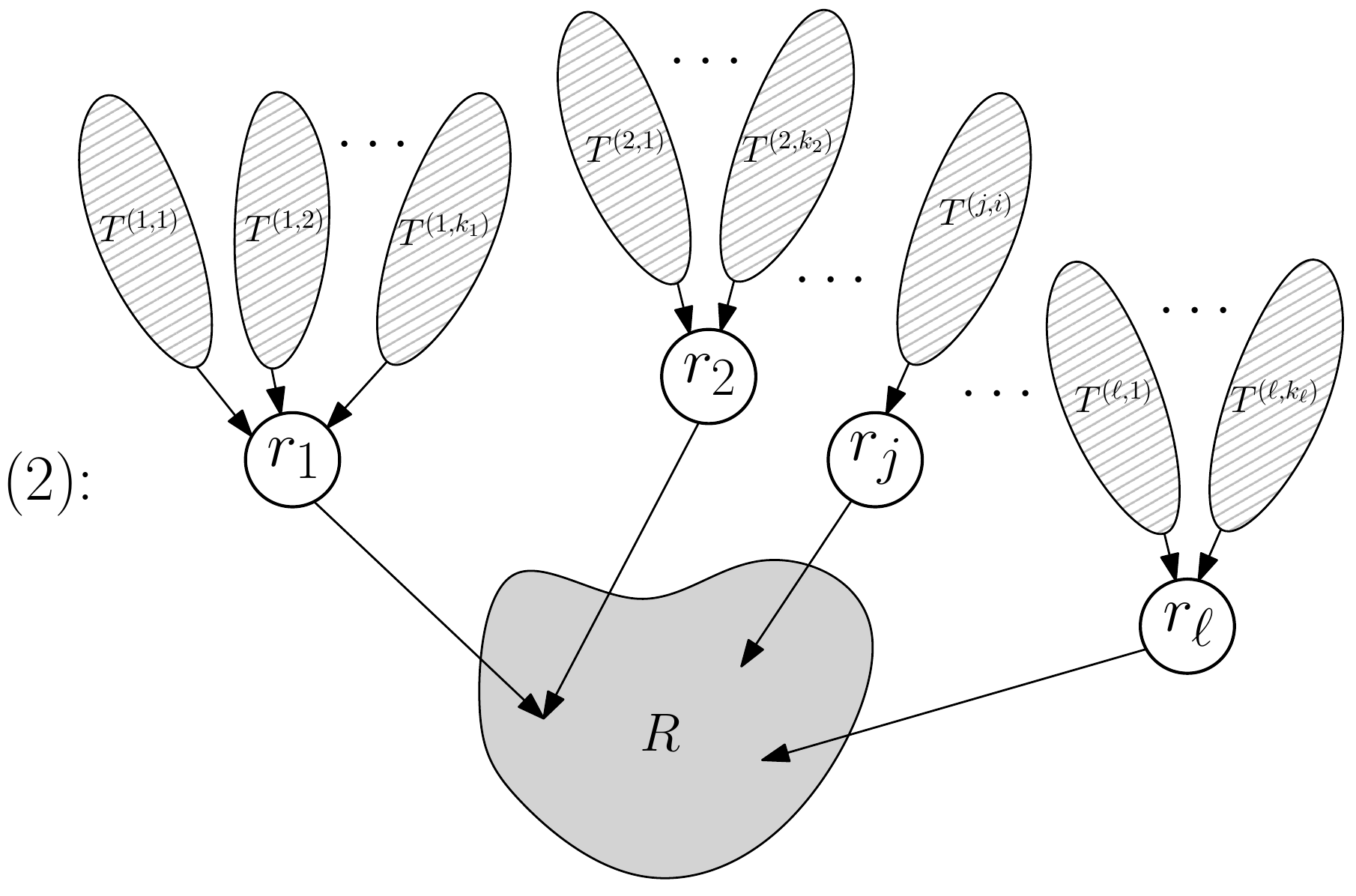}
\caption{Decomposition of a general tree parking function $(T,s)$ w.r.t.\ either the unoccupied root $r$ or the root-cluster $R$.\label{pic:generaltreedecomposition}}
\end{center}
\end{figure}

We consider now the general case of $(n,m)$-tree parking functions and $(n,m)$-tree parking distributions for size-$n$ trees and sequences (resp.\ multisets) of $m \le n$ drivers. Our enumerative approach relies on a combinatorial decomposition of a $(n,m)$-tree parking function $(T,s)$ w.r.t.\ either the root of $T$ (if it is unoccupied) or the maximal subtree $R$ (let us call it the ``root-cluster'') of $T$ that contains the root of $T$ and that contains only nodes occupied during the parking procedure. Depending on whether the root $r$ of $T$ is occupied during the parking procedure or not, we have to distinguish between two cases:
\begin{enumerate}
\item Root $r$ is not occupied during the parking procedure: decomposing $T$ into $r$ and the subtrees $T^{(1)}, \dots, T^{(k)}$ linked to $r$, and denoting by $s^{(1)}, \dots, s^{(k)}$ the corresponding subsequences of $s$ of drivers arriving (and parking) in the respective subtrees, it holds that $(T^{(j)},s^{(j)})$, $1 \le j \le k$, are themselves general tree parking functions of size $|T^{(j)}|$ and length $|s^{(j)}|$, $1 \le j \le k$.
\item Root $r$ is occupied during the parking procedure: removing the root-cluster $R$, $T$ decomposes into trees $T^{(1)}, \dots, T^{(\ell)}$. For each of these trees $T^{(j)}$, $1 \le j \le \ell$, let us denote by $r_{j}$ its root and by $T^{(j,1)}, \dots, T^{(j,k_{j})}$ the subtrees linked to $r_{j}$; then it holds that each $r_{j}$ is unoccupied during the parking procedure. When we denote by $s^{(R)}$ and by $s^{(j,i)}$, $1 \le j \le \ell$ and $1 \le i \le k_{j}$, the subsequences of $s$ of drivers arriving (and parking) in the root-cluster $R$ or the corresponding subtrees $T^{(j,i)}$, respectively, it holds that $(R,s^{(R)})$ is a tree parking function of size $|R|$, whereas each $(T^{(j,i)},s^{(j,i)})$ is a general tree parking function of size $|T^{(j,i)}|$ and length $|s^{(j,i)}|$.
\end{enumerate}
The decomposition is illustrated in Figure~\ref{pic:generaltreedecomposition} and holds in a completely analogous way also for general tree parking distributions.

\subsection{Ordered trees\label{ssec:GeneralParking_Ordered}}

\subsubsection{General tree parking functions}
Here we treat $(n,m)$-tree parking functions and distributions for ordered trees and start with parking functions. Thus let us denote by $\mathcal{F}$ the family of general tree parking functions and by $\mathcal{G}$ the family of parking functions for ordered trees. We will assume that for each general tree parking function $f \in \mathcal{F}$ every unoccupied node in the tree component is marked with a marker $U$, and that for each parking function $g \in \mathcal{G}$ every attachment point in the tree component is marked with a marker $A$. We may express via $\mathcal{G} = \mathcal{G}(\mathcal{Z},A)$ the possession of nodes $\mathcal{Z}$ and attachment points $A$ for each object of $\mathcal{G}$. The decomposition of a general tree parking function w.r.t.\ the root-cluster as given in Section~\ref{ssec:GeneralParking_Decomposition} and illustrated in Figure~\ref{pic:generaltreedecomposition} leads to the formal equation
\begin{equation}
  \mathcal{F} = \mathcal{Z} \times \{U\} \ast \textsc{Seq}(\mathcal{F}) + \mathcal{G}\big(\mathcal{Z}, \textsc{Seq}\big(\mathcal{Z} \times \{U\} \ast \textsc{Seq}(\mathcal{F})\big)\big),
	\label{eqn:general_parking_formal_G-F}
\end{equation}
where the two cases of the decomposition reflect the corresponding summands in this formal equation. For the second summand we use that at each attachment point in the root-cluster there can be attached a sequence of trees, and each such tree consists of an unoccupied root node to which a sequence of subtrees (corresponding to the tree components of respective general tree parking functions) is attached.

We introduce the (double-exponential) bivariate generating function $F(z,u)$, where the variable $u$ counts the number of unoccupied nodes,
\begin{equation*}
  F:=F(z,u) = \sum_{f=(T,s) \in \mathcal{F}} \frac{z^{|T|} u^{|T|-|s|}}{|T|! \, |s|!} = \sum_{n \ge 1} \sum_{0 \le m \le n} F_{n,m} \frac{z^{n} u^{n-m}}{n! \, m!},
\end{equation*}
with $F_{n,m} = |\{f=(T,s) \in \mathcal{F} : |T|=n \; \text{and} \; |s|=m\}|$ the total number of $(n,m)$-tree parking functions for ordered trees. Furthermore we require the generating function $G(z)$ of the number of tree parking functions for ordered trees as defined in \eqref{eqn:ordered_def_G}. By taking into account that any ordered tree $T$ of size $n$ has exactly $2n-1$ attachment points, the formal equation~\eqref{eqn:general_parking_formal_G-F} leads at the level of generating functions to the following relation between $F = F(z,u)$ and $G = G(z)$:
\begin{equation*}
  F = \frac{zu}{1-F} + \Big(1-\frac{zu}{1-F}\Big) \cdot G\bigg(\frac{z}{\big(1-\frac{zu}{1-F}\big)^{2}}\bigg),
\end{equation*}
resp.\
\begin{equation}
  \frac{F - \frac{zu}{1-F}}{1-\frac{zu}{1-F}} = G\bigg(\frac{z}{\big(1-\frac{zu}{1-F}\big)^{2}}\bigg).
	\label{eqn:general_parking_eqn_G-F}
\end{equation}
It is slightly tricky to obtain from \eqref{eqn:general_parking_eqn_G-F} a suitable characterization of the function $F(z,u)$. We start with the characterization \eqref{eqn:ordered_feq_G-Q} of the solution of $G=G(z)$ via an auxiliary function $Q := Q(z)$ and evaluate them at $\frac{z}{\big(1-\frac{zu}{1-F}\big)^{2}}$. Introducing
$\hat{G}:=G\Big(\frac{z}{\big(1-\frac{zu}{1-F}\big)^{2}}\Big)$ and $\hat{Q}:=Q\Big(\frac{z}{\big(1-\frac{zu}{1-F}\big)^{2}}\Big)$,
yields
\begin{equation}
  \hat{G} = 1 - \frac{1}{(1+\hat{Q}) (1-\hat{Q}) e^{\hat{Q}}} \quad \text{and} \quad \hat{Q} = \frac{z (1+\hat{Q})^{2} e^{\hat{Q}}}{(1-\frac{zu}{1-F})^{2}}.
	\label{eqn:general_parking_rel_G-Q}
\end{equation}
Thus, by using \eqref{eqn:general_parking_eqn_G-F} and \eqref{eqn:general_parking_rel_G-Q}, we obtain
$\frac{F - \frac{zu}{1-F}}{1-\frac{zu}{1-F}} = \frac{F-1}{1-\frac{zu}{1-F}} +1 = \hat{G} = 1- \frac{1}{(1+\hat{Q}) (1-\hat{Q}) e^{\hat{Q}}}$, and further
\begin{equation}
  \frac{1-F}{1-\frac{zu}{1-F}} = \frac{1}{(1+\hat{Q}) (1-\hat{Q}) e^{\hat{Q}}}.
	\label{eqn:general_parking_rel1_F-Q}
\end{equation}
Taking squares and rearranging yields $\frac{1}{(1-\frac{zu}{1-F})^{2}} = \frac{1}{(1-F)^{2} (1+\hat{Q})^{2} (1-\hat{Q})^{2} \, e^{2\hat{Q}}}$; plugging this into the relation \eqref{eqn:general_parking_rel_G-Q} for $\hat{Q}$ gives
\begin{equation}
  \hat{Q} = \frac{z}{(1-F)^{2} (1-\hat{Q})^{2} e^{\hat{Q}}},
	\label{eqn:general_parking_rel2_F-Q}
\end{equation}
and further
\begin{equation}
  \frac{zu}{(1-F)^{2}} = u \hat{Q} (1-\hat{Q})^{2} e^{\hat{Q}}.
	\label{eqn:general_parking_rel3_F-Q}
\end{equation}
Rearranging \eqref{eqn:general_parking_rel1_F-Q} to $\frac{1}{1-F} - \frac{zu}{(1-F)^{2}} = (1+\hat{Q}) (1-\hat{Q}) e^{\hat{Q}}$ and plugging \eqref{eqn:general_parking_rel3_F-Q} into it, leads to the relation
\begin{equation}
  \frac{1}{1-F} = (1-\hat{Q}) \big(1+\hat{Q}+u\hat{Q}(1-\hat{Q})\big) e^{\hat{Q}},
	\label{eqn:general_parking_rel4_F-Q}
\end{equation}
which easily gives the desired solution for $F = F(z,u)$ via the auxiliary function $\hat{Q} = \hat{Q}(z,u)$. Namely, plugging \eqref{eqn:general_parking_rel4_F-Q} into \eqref{eqn:general_parking_rel2_F-Q} shows that $\hat{Q}$ is characterized implicitly via the functional equation
\begin{equation}
  \hat{Q} = z \big(1+\hat{Q} + u\hat{Q}(1-\hat{Q})\big)^{2} e^{\hat{Q}},
\end{equation}
and according to \eqref{eqn:general_parking_rel4_F-Q} $F$ is related to it by means of
\begin{equation}
  F = 1 - \frac{1}{(1-\hat{Q}) \big(1+\hat{Q}+u\hat{Q}(1-\hat{Q})\big) e^{\hat{Q}}},
\end{equation}
which is the solution stated in Theorem~\ref{thm:gf_general_treeparking}.

\subsubsection{General tree parking distributions}
These considerations can be applied also to gain results for general tree parking distributions. Let $\tilde{\mathcal{F}}$ and $\tilde{\mathcal{G}}$ denote the family of general tree parking distributions and tree parking distributions, respectively, for ordered trees, where we express via $\tilde{\mathcal{G}} = \tilde{\mathcal{G}}(\mathcal{Z},A)$ the dependence of any object of $\tilde{\mathcal{G}}$ by its nodes $\mathcal{Z}$ and attachment points $A$. The decomposition of a general tree parking distribution w.r.t.\ the root-cluster leads to an analogue of the formal equation \eqref{eqn:general_parking_formal_G-F} relating the families $\tilde{\mathcal{F}}$ and $\tilde{\mathcal{G}}$.

Introducing the bivariate generating function
\begin{equation*}
  \tilde{F} := \tilde{F}(z,u) = \sum_{\tilde{f}=(T,\tilde{s}) \in \tilde{\mathcal{F}}} \frac{z^{|T|} u^{|T| - |\tilde{s}|}}{|T|!} = \sum_{n \ge 1} \sum_{0 \le m \le n} \tilde{F}_{n,m} \frac{z^{n} u^{n-m}}{n!},
\end{equation*}
with $\tilde{F}_{n,m} = |\{\tilde{f}=(T,\tilde{s}) \in \tilde{\mathcal{F}} : |T|=n \; \text{and} \; |\tilde{s}|=m\}|$ the total number of $(n,m)$-tree parking distributions for ordered trees, and using the generating function $\tilde{G}(z)$ of the number of tree parking distributions as defined in \eqref{eqn:ordered_distribution_def_G}, we obtain from the formal equation~\eqref{eqn:general_parking_formal_G-F} the following relation between $\tilde{F} = \tilde{F}(z,u)$ and $\tilde{G} = \tilde{G}(z)$:
\begin{equation}
  \frac{\tilde{F} - \frac{zu}{1-\tilde{F}}}{1-\frac{zu}{1-\tilde{F}}} = \tilde{G}\bigg(\frac{z}{\big(1-\frac{zu}{1-\tilde{F}}\big)^{2}}\bigg).
	\label{eqn:general_distribution_eqn_F-G}
\end{equation}

Using the characterization \eqref{eqn:ordered_distribution_sol_G_Q} of the generating function $\tilde{G}$ via an auxiliary function $Q$ and
introducing $\hat{G} := \tilde{G}\Big(\frac{z}{\big(1-\frac{zu}{1-\tilde{F}}\big)^{2}}\Big)$ and $\hat{Q} := Q\Big(\frac{z}{\big(1-\frac{zu}{1-\tilde{F}}\big)^{2}}\Big)$, yields
\begin{equation}
  \hat{G} = 1 - \frac{1}{(1+\hat{Q})^{2} \, (1-\hat{Q})} \quad \text{and} \quad \hat{Q} = \frac{z (1+\hat{Q})^{4}}{(1-\frac{zu}{1-\tilde{F}})^{2}}.
	\label{eqn:general_distribution_rel_G-Q}
\end{equation}
Using \eqref{eqn:general_distribution_rel_G-Q}, we get from \eqref{eqn:general_distribution_eqn_F-G} the relation
\begin{equation*}
  \frac{1-\tilde{F}}{1-\frac{zu}{1-\tilde{F}}} = \frac{1}{(1+\hat{Q})^{2} (1-\hat{Q})},
\end{equation*}
and after a few further steps (similar to the ones carried out for tree parking functions) show the result given in Theorem~\ref{thm:gf_general_treeparking}, namely that the solution of $\tilde{F} = \tilde{F}(z,u)$ can be given by means of the auxiliary function $\hat{Q} = \hat{Q}(z,u)$, which is itself characterized implicitly via a functional equation:
\begin{equation}
  \tilde{F} = 1 - \frac{1}{(1-\hat{Q}) \big((1+\hat{Q})^{2} + u\hat{Q}(1-\hat{Q})\big)}, \qquad \hat{Q} = z \big((1+\hat{Q})^{2} +u\hat{Q}(1-\hat{Q})\big)^{2}.
\end{equation}

\subsection{Combinatorial tree families}

The way of getting enumerative results for general tree parking functions and distributions from the corresponding solutions for tree parking functions and distributions, respectively, carried out for ordered trees in Section~\ref{ssec:GeneralParking_Ordered} can be extended to the other combinatorial tree families. Again we will only sketch them for parking functions, since the computations for parking distributions can be carried out in an analogous way.

The decomposition of a general tree parking function w.r.t.\ the root-cluster can be translated into a formal equation relating the family of general tree parking functions $\mathcal{F}$ and the family of tree parking functions $\mathcal{G} = \mathcal{G}(\mathcal{Z},A)$. In case that the root node is occupied one has to take into account that for unordered trees at each attachment point in the root-cluster there is attached a (possibly empty) set of general tree parking functions, whereas for $d$-ary trees at each attachment point of the root-cluster there is either attached a general tree parking function or not; for $d$-bundled trees, as for the particular instance of ordered trees, at each attachment point of the root-cluster there is attached a sequence of general tree parking functions. This leads to the following symbolic equations:
\begin{align}
\text{unordered trees:} & \quad \mathcal{F} = \mathcal{Z} \times \{U\} \ast \textsc{Set}(\mathcal{F}) + \mathcal{G}\big(\mathcal{Z},\textsc{Set}\big(\mathcal{Z} \times \{U\} \ast \textsc{Set}(\mathcal{F})\big)\big),\notag\\
\text{$d$-ary trees:} & \quad \mathcal{F} = \mathcal{Z} \times \{U\} \ast (\{\epsilon\}+\mathcal{F})^{d} + \mathcal{G}\big(\mathcal{Z},\{\epsilon\}+\mathcal{Z}\times \{U\} \ast (\{\epsilon\}+\mathcal{F})^{d}\big),\label{eqn:general_parking_combinatorial_formal_G-F}\\
\text{$d$-bundled trees:} & \quad \mathcal{F} = \mathcal{Z} \times \{U\} \ast \big(\textsc{Seq}(\mathcal{F})\big)^{d} + \mathcal{G}\big(\mathcal{Z},\textsc{Seq}\big(\mathcal{Z} \times \{U\} \ast (\textsc{Seq}(\mathcal{F}))^{d}\big)\big).\notag
\end{align}

Introducing the generating functions
\begin{equation*}
  F := F(z,u) = \sum_{f=(T,s) \in \mathcal{F}} \frac{z^{|T|} u^{|T| - |s|}}{|T|! \, |s|!} = \sum_{n \ge 1} \sum_{0 \le m \le n} F_{n,m} \frac{z^{n} u^{n-m}}{n! \, m!},
\end{equation*}
with $F_{n,m} = |\{f=(T,s) \in \mathcal{F} : |T|=n \; \text{and} \; |s|=m\}|$ the total number of $(n,m)$-tree parking functions for the corresponding tree family (and keeping in mind the number of attachment points of a size-$n$ tree of the respective tree families), we get from \eqref{eqn:general_parking_combinatorial_formal_G-F} the following relations between the generating functions $G = G(z)$ (for tree parking functions as defined in Section~\ref{ssec:Parking_Function_Combinatorial}) and $F = F(z,u)$:
\begin{align}
\text{unordered trees:} & \quad F = zu e^{F} + G\Big(z e^{zue^{F}}\Big),\notag\\
\text{$d$-ary trees:} & \quad F = zu (1+F)^{d} + \big(1+zu(1+F)^{d}\big) \cdot G\left(z \big(1+zu(1+F)^{d}\big)^{d-1}\right),\label{eqn:general_parking_combinatorial_eqn_G-F}\\
\text{$d$-bundled trees:} & \quad F = \frac{zu}{(1-F)^{d}} + \Big(1-\frac{zu}{(1-F)^{d}}\Big) \cdot G\bigg(\frac{z}{(1-\frac{zu}{(1-F)^{d}})^{d+1}}\bigg).\notag
\end{align}
Computations analogous to the ones for ordered trees, by taking into account the characterizations of the g.f.\ for $G(z)$ collected in Table~\ref{tab:gf_treeparking}, show the results stated in Theorem~\ref{thm:gf_general_treeparking}. Again it is not difficult to show that the functions given there satisfy above equations \eqref{eqn:general_parking_combinatorial_eqn_G-F} and are thus indeed the required solutions.

\section{Asymptotic results\label{sec:Asymptotic}}

\subsection{Univariate sequences}

The asymptotic behaviour of the sequences $G_{n}$, $P_{n}$, $\tilde{G}_{n}$, and $\tilde{P}_{n}$, for $n \to \infty$, enumerating different kinds of parking functions for combinatorial tree families stated in Theorem~\ref{thm:asymptotic_treeparking} resp.\ Table~\ref{tab:asymptotic_treeparking} can be obtained from the generating functions solutions given in Theorem~\ref{thm:gf_treeparking} via a standard application of singularity analysis to these implicitly defined functions. Since the auxiliary functions $Q = Q(z)$ are all characterized by a functional equation $Q = z \varphi(Q)$ satisfying the conditions for so-called singular inversion for aperiodic functions $\varphi$ as formulated in \cite[Theorem~VI.6]{FlaSed2009}, we obtain directly that $Q(z)$ has a unique dominant singularity, i.e., singularity of smallest modulus, $\rho = \frac{\tau}{\varphi(\tau)}$, with $\tau$ given by the smallest positive real root of the equation $\varphi(\tau) = \tau \varphi'(\tau)$, and that $Q(z)$ admits a singular expansion around $\rho$ of the form
\begin{equation}
  Q = \tau - {\textstyle{\sqrt{\frac{2 \varphi(\tau)}{\varphi''(\tau)}} \sqrt{1-\frac{z}{\rho}}}} + \sum_{k \ge 2} {\textstyle{d_{k} \big(1-\frac{z}{\rho}\big)^{\frac{k}{2}}}},
	\label{eqn:asymptotic_univariate_expansion_Q}
\end{equation}
with computable constants $d_{k}$. The g.f.\ $X := X(z)$ (meaning $G(z)$, $P(z)$, $\tilde{G}(z)$, or $\tilde{P}(z)$) of interest can be expressed via these auxiliary function $Q(z)$, and one easily obtains that $X(z)$ and the corresponding $Q(z)$ have the same unique dominant singularity $\rho$ and that the singular expansion of $Q$ easily yields a singular expansion of $X$. There is the interesting aspect that in all cases considered the coefficient of the square root-term vanishes, i.e., that we get an expansion of the form 
\begin{equation*}
  X = c_{0} + {\textstyle{c_{2} \, \big(1-\frac{z}{\rho}\big)}} + \sum_{k \ge 3} {\textstyle{c_{k} \big(1-\frac{z}{\rho}\big)^{\frac{k}{2}}}}.
\end{equation*}
An application of the transfer theorems of Flajolet and Odlyzko (see \cite{FlaSed2009}) to this expansion show then an asymptotic behaviour of the coefficients $[z^{n}] X(z) \sim C n^{-\frac{5}{2}} \rho^{-n}$, with $C = \frac{3 c_{3}}{4 \sqrt{\pi}}$. However, due to these cancellations in the singular expansions of $X$, which would require to determine also the coefficients $d_{2}$ and $d_{3}$ in the expansion \eqref{eqn:asymptotic_univariate_expansion_Q} of $Q$, it seems to be easier in practice to expand the derivative $X'(z)$ of the g.f.\ of interest around $\rho$; then the first two terms in the expansion \eqref{eqn:asymptotic_univariate_expansion_Q} are sufficient.

\smallskip

We exemplify the computations for the particular instance of tree parking functions for ordered trees. As stated in Theorem~\ref{thm:gf_treeparking}, the auxiliary function $Q = Q(z)$ satisfies the functional equation $Q= z \varphi(Q)$, with $\varphi(Q) = (1+Q)^{2} e^{Q}$. According to the implicit function theorem, the equation $h(z,Q) := Q - z \varphi(Q) = 0$ cannot be resolved w.r.t.\ $Q$ in a locally unique way for points $(z,Q)=(\rho,\tau)$ satisfying $h(\rho,\tau)=0$ and $h_{Q}(\rho,\tau)=0$, which leads to the given characterization $\varphi(\tau) = \tau \varphi'(\tau)$, and $\rho = \frac{\tau}{\varphi(\tau)} = \frac{1}{\varphi'(\tau)}$. In the present case we get $(1+\tau)^{2} e^{\tau} = \tau (1+\tau)(3+\tau)e^{\tau}$ with the only positive real root $\tau = \sqrt{2}-1$ (the restriction to the positive real root is justified, since the coefficients $Q_{n}$ of the power series $Q(z) = \sum_{n \ge 1} Q_{n} z^{n}$ are all non-negative), which implies $\rho = \frac{\sqrt{2}-1}{2 e^{\sqrt{2}-1}}$. For the local expansion of $Q$ around $\rho$ we compute $\frac{2 \varphi(\tau)}{\varphi''(\tau)} = \tau$, thus get
\begin{equation}
  Q = {\textstyle{\tau - \sqrt{\tau} \sqrt{1-\frac{z}{\rho}} + \mathcal{O}\big(1-\frac{z}{\rho}\big)}}.
	\label{eqn:asymptotic_univariate_ordered_expansion_Q}
\end{equation}
From Theorem~\ref{thm:gf_treeparking} we obtain that $G = G(z)$ is given by $G=1-\frac{1}{(1+Q)(1-Q)e^{Q}}$. Simple computations show that $G'(z)$ can be expressed by means of $Q$ via $G'(z) = \frac{1+Q}{(1-Q)^{2}}$, thus $G'$ also has its unique dominant singularity at $z=\rho$. Plugging above local expansion \eqref{eqn:asymptotic_univariate_ordered_expansion_Q} of $Q$ into this representation of $G'$ yields after some computations the local expansion
\begin{equation*}
  G' = {\textstyle{2+\frac{3 \sqrt{2}}{2} - \sqrt{\sqrt{2}-1} \, \big(\frac{17}{2} + 6 \sqrt{2}\big) \sqrt{1-\frac{z}{\rho}} + \mathcal{O}\big(1-\frac{z}{\rho}\big)}},
\end{equation*}
and an application of transfer theorems gives
\begin{equation*}
  [z^{n}] G(z) = \frac{1}{n} [z^{n-1}] G'(z) \sim \frac{\sqrt{\sqrt{2}-1} \, \big(\frac{17}{2}+6\sqrt{2}\big)}{2 \sqrt{\pi}} \, n^{-\frac{5}{2}} \rho^{-n+1},
\end{equation*}
and furthermore the asymptotic result stated in Theorem~\ref{thm:asymptotic_treeparking}:
\begin{equation*}
  G_{n} = (n!)^{2} [z^{n}] G(z) \sim (n!)^{2} \cdot C \, n^{-\frac{5}{2}} \rho^{-n}, \quad \text{with} \quad C = \frac{(\sqrt{2}-1)^{\frac{3}{2}} (17+12\sqrt{2})}{8 \, e^{\sqrt{2}-1} \sqrt{\pi}}.
\end{equation*}

\subsection{Bivariate sequences}

From the representations of the generating functions $F(z,u)$ and $\tilde{F}(z,u)$, respectively, given in Theorem~\ref{thm:gf_general_treeparking}, the asymptotic behaviour of the coefficients can be obtained via saddle point techniques. For the particular instance of tree parking functions for unordered trees this approach has been carried out in detail in \cite{LackPan2016} and it can be adapted easily to the other tree families. Thus in the following we only give a raw sketch of the proof of the result stated in Theorem~\ref{thm:asymptotic_general_treeparking} by exemplifying the computations for the instance of tree parking distributions for ordered trees.

Namely, from the representation $\tilde{F} = 1-\frac{1}{(1-Q)((1+Q)^{2}+uQ(1-Q))}$, with $Q=z \varphi(Q,u) = z \big((1+Q)^{2}+uQ(1-Q)\big)^{2}$, as stated in Table~\ref{tab:gf_general_treeparking} (where for a better readability we write here $Q$ instead of $\hat{Q}$), thus $\frac{d\tilde{F}}{dQ} = \frac{(1-3Q)(1+Q+u(1-Q))}{(1-Q)^{2} ((1+Q)^{2}+uQ(1-Q))^{2}}$, one gets, by an application of the Lagrange-B{\"u}rmann inversion formula, after straightforward computations the following representation of the number of parking distributions $\tilde{F}_{n,m}$:
\begin{align*}
  \tilde{F}_{n,m} & = n! [z^{n} u^{n-m}] \tilde{F}(z,u) = (n-1)! [Q^{n-1} u^{n-m}] {\textstyle{\frac{d\tilde{F}}{dQ}}} \cdot \big(\varphi(Q,u)\big)^{n}\\
	& = (n-1)! \binom{2n-1}{n-m} [Q^{m}] (1-3Q) \Big(\frac{n-m}{2n-1} + Q\Big) (1-Q)^{n-m-2} (1+Q)^{2n+2m-3}.
\end{align*}
From this representation one could easily obtain (as it also holds for other tree families) an explicit formula for the numbers $\tilde{F}_{n,m}$, but we omit it here and continue with asymptotic considerations. To this aim we introduce the auxiliary sequence $A_{n,m}$ via
\begin{equation*}
  A_{n,m} := [Q^{m}] (1-3Q) \Big(\frac{n-m}{2n-1} + Q\Big) (1-Q)^{n-m-2} (1+Q)^{2n+2m-3}
	= \frac{1}{2 \pi i} \oint g(w) \cdot e^{n h(w)} dw,
\end{equation*}
with
\begin{equation*}
  g(w) = \frac{(1-3w)(\frac{n-m}{2n-1}+w)}{w(1-w)(1+w)^{3}}, \quad h(w) = \big(1-\frac{m}{n}\big)\ln(1-w) +2\big(1+\frac{m}{n}\big)\ln(1+w)-\frac{m}{n}\ln w,
\end{equation*}
and where we choose as contour in the contour integral representation a suitable simple positively oriented closed curve around the origin. To locate the saddle points in the relevant part $e^{n h(w)}$ of the integrand, we have to solve the equation 
\begin{equation*}
  h'(w) = \frac{(3w-1)(m-nw)}{nw(1-w)(1+w)} = 0,
\end{equation*}
yielding the solutions $w_{1}=\frac{m}{n}$ and $w_{2}=\frac{1}{3}$. The asymptotic behaviour of $A_{n,m}$ (and thus of $\tilde{F}_{n,m}$) depends on which of the two saddle points is the smaller one (or whether they coalesce). We only treat here the case $w_{1} < w_{2}$ (for the other cases we refer to \cite{LackPan2016}) and assume that the load factor satisfies $\alpha = \frac{m}{n} \le \alpha_{0} -\epsilon = \frac{1}{3} - \epsilon$, for an arbitrary small but fixed $\epsilon > 0$.
Thus, we choose as contour a circle centered at the origin with radius $\alpha$, and thus passing through the saddle point $w_{1}$. Using the parametrization $\Gamma = \{w = \alpha e^{i t} : t \in [-\pi, \pi]\}$ for the curve in the integral representation of $A_{n,m}$, we get
\begin{equation*}
  A_{n,m} = \frac{1}{2\pi} \int_{-\pi}^{\pi} \alpha e^{it} g(\alpha e^{it}) \, e^{n h(\alpha e^{it})} dt.
\end{equation*}
With the local expansions around $t=0$:
\begin{align*}
  \alpha e^{it} g(\alpha e^{it}) & = \frac{1-3\alpha}{2(1+\alpha)^{2} (1-\alpha)^{2}} \cdot {\textstyle{\big(1+\mathcal{O}(\frac{m}{n}t)\big)}},\\
	e^{n h(\alpha e^{it})} & = \frac{(1-\alpha)^{n-m} (1+\alpha)^{2(n+m)}}{\alpha^{m}} \, e^{-\frac{m(1-3\alpha)}{2(1-\alpha)(1+\alpha)}t^{2}} \cdot \big(1+\mathcal{O}(mt^{3})\big),
\end{align*}
one eventually obtains the asymptotic evaluation
\begin{align*}
  A_{n,m} & \sim \frac{(1-\alpha)^{n-m} (1+\alpha)^{2(n+m)}}{\alpha^{m}} \cdot \frac{1-3\alpha}{2 (1-\alpha)^{2} (1+\alpha)^{2}} \cdot \frac{1}{2\pi} \int_{-\infty}^{+\infty} e^{-\frac{m(1-3\alpha)}{2(1-\alpha)(1+\alpha)}t^{2}} dt\\
	& \sim \frac{(1-\alpha)^{n-m} (1+\alpha)^{2(n+m)}}{\alpha^{m}} \cdot \frac{\sqrt{2} \sqrt{1-3\alpha}}{4 \sqrt{\pi} \, (1-\alpha)^{\frac{3}{2}} (1+\alpha)^{\frac{3}{2}} \sqrt{\alpha} \sqrt{n}}.
\end{align*}
We are actually interested in the probability $p_{n,m} = \frac{\tilde{F}_{n,m}}{T_{n} \binom{n+m-1}{m}} = \frac{(n-1)! \binom{2n-1}{n-m} A_{n,m}}{(n-1)! \binom{2n-2}{n-1} \binom{n+m-1}{m}}$ that a randomly chosen pair $(T,\tilde{s})$ with $T$ a size-$n$ ordered tree and $\tilde{s}$ a multiset of size $m$ on $[n]$ is a $(n,m)$-tree parking distribution, for which we obtain after a standard application of Stirling's formula for the factorials (see, e.g., \cite{FlaSed2009}) the asymptotic evaluation stated in Table~\ref{tab:asymptotic_general_treeparking}:
\begin{equation*}
  p_{n,m} \sim \frac{\sqrt{1-3\alpha}}{(1-\alpha)^{2} \sqrt{1+\alpha}}, \quad \text{for \, ${\textstyle{0 \le \alpha \le \frac{1}{3} - \epsilon}}$}.
\end{equation*}

\section{Outlook and open problems\label{sec:Outlook}}

As shown, the generating functions approach presented based on combinatorial decompositions is capable of treating the enumeration problem for various kinds of tree parking functions and tree families. It could be extended further, where we mention a few such directions. 

Although we have restricted in this work mainly on what we have called combinatorial tree families, the approach also works for further tree families, and we have given a few such examples in Section~\ref{ssec:parking_distribution_further_trees} on tree parking distributions. However, it would be of particular interest to treat also tree families that do not fall into the class of simple families of trees. One such family are so-called recursive trees, i.e., labelled unordered trees, where the labels along each leaf-to-root path are increasing. For this tree family the approach presented could be applied, but yields a second-order non-linear differential equation, which does not seem to be easily amenable.

Another promising direction is to study defective tree parking functions, i.e., to study the ``overflow'' of the number of drivers, which could not park successfully (see \cite{CamJohPreSchwe2008} for a respective treatment on ordinary parking functions). With the methods proposed it seems that the results presented could be extended to this problem and the author plans to comment on that in a future work.

Further possible directions for which, at least in principle, the approach presented could be applied, concern, e.g., restrictions on the number of arrivals on parking spaces (e.g., to at most two arrivals or to an even number of arrivals), ``friends parking'' (pairs, or more general $j$-tuples, of drivers always arrive at the same preferred parking space), or ``truck parking'' (drivers require two, or more general $j$, consecutive parking spaces).

Of course, there are further problems in connection with parking on trees, which do not seem to be amenable with the presented approach. One interesting such problem is to consider parking functions on trees, where the edges are oriented towards the leaves. Although the paths of the drivers in search of a free parking space are not uniquely determined, the notion makes sense as has been shown in \cite{KinYan2020}. Treating the enumeration problem for such kinds of parking functions is an open problem.

As mentioned earlier, the asymptotic behaviour of the enumeration sequences for parking function varieties reminds one on the asymptotic behaviour of some enumeration problems for maps. Since there is even an enumeration sequence (prime tree parking distributions on ordered trees), which occurs in the enumeration of a certain kind of maps, it would be interesting to establish links between these combinatorial structures.

\end{document}